%% file: DiscForms.tex
\documentclass[a4paper,12pt,oneside]{article}
\usepackage[latin1]{inputenc}
\usepackage[sumlimits,namelimits]{amsmath}
\usepackage{amssymb,amscd,latexsym}
\usepackage{theorem}
\usepackage{enumerate}
\usepackage[T1]{fontenc}
\usepackage{graphicx}
\usepackage{multicol}
\usepackage{pstcol,pst-char}

\newtheorem{dfn}{Definition}[section]

\newtheorem{thm}[dfn]{Theorem}

\newtheorem{lmm}[dfn]{Lemma}
\newtheorem{rmk}[dfn]{Remark}

\def\eps{\varepsilon}
\def\R{\mathbb{R}}
\def\N{\mathbb{N}}
\def\U{\mathcal{U}}
\def\D{\mathcal{D}}
\def\S{\mathcal{S}}

\def\Cech{\v{C}ech }
\def\K{\mathbb{M}(n,\kappa,r_0)}

\newcommand{\CqR}[2]{\mathcal{C}^{#1}(#2)}
\newcommand{\Cq}[3]{\mathcal{C}^{#1}(\mathcal{#2},#3)}

\numberwithin{equation}{section}

\title{Discretization of Riemannian manifolds applied to the Hodge Laplacian}
\author{Tatiana Mantuano
\thanks{Supported by Swiss National Science Foundation, grant No. 20-109130.}
}
\date{}

\begin{document}

\maketitle

{\abstract For $\kappa \geq 0$ and $r_0 > 0$, let $\K$ be the set
of all connected compact $n$-dimensional Riemannian manifolds such
that $|K_g| \leq \kappa$ and $Inj(M,g) \geq r_0$. We study the
relation between the $k^{\text{th}}$ positive eigenvalue of the
Hodge Laplacian on differential forms and the $k^{\text{th}}$
positive eigenvalue of the combinatorial Laplacian associated to
an open cover (acting on \Cech cochains). We show that for a fixed
sufficiently small $\eps > 0$ there exist positive constants $c_1$
and $c_2$ depending only on $n$, $\kappa$ and $\eps$ such that for
any $M \in \K$ and for any $\eps$-discretization $X$ of $M$ we
have $c_1 \lambda_{k,p}(X) \leq \lambda_{k,p}(M) \leq c_2
\lambda_{k,p}(X)$ for any $k \leq K$ ($K$ depends on $X$).
Moreover, we find a lower bound for the spectrum of the
combinatorial Laplacian and a lower bound for the spectrum of the
Hodge Laplacian.

\smallskip


\smallskip


}


\input{DiscFormsIntro}

\input{DiscFormsSettings}

\input{DiscFormsMain}

\input{DiscFormsApplication}

\input{DiscFormsAppendix}

\input{DiscFormsReferences}

\end{document}

%% file: DiscFormsIntro.tex
\section{Introduction}

Several works like \cite{Br1}, \cite{Bur}, \cite{Buser} and more
recently \cite{Ma} have shown that discretizing a Riemannian
manifold may be really powerful in order to study the spectrum of
the Laplacian acting on functions. The question we want to answer
here is "Is there a similar tool for understanding the spectrum of
the Hodge Laplacian ($\Delta = dd^* + d^*d$) acting on
differential forms?". Part of an answer is given by the de Rham
Theorem (saying that the de Rham cohomology of a compact manifold
is isomorphic to the singular cohomology and to the \Cech
cohomology) and several authors have been more or less inspired by
this theorem to study the spectrum of $\Delta$. For instance, in
\cite{DP}, Dodziuk and Patodi show that for a fixed compact
Riemannian manifold, we can approximate the spectrum of the Hodge
Laplacian with the spectrum of a combinatorial Laplacian
associated to finer and finer triangulations of the manifold. The
main idea in their proof is to associate \Cech cochains to smooth
forms and vice versa via the integration on simplices and via the
Whitney map. Both tools are really crucial in the proof of the de
Rham Theorem as they induce the isomorphism between de Rham
cohomology and singular cohomology. In \cite{CT} and in \cite{MG},
the authors use another proof of de Rham Theorem due to A. Weil
and based on the \Cech - de Rham double complexe (see \cite{Go}).
In \cite{CT}, Chanillo and Trèves bound from below the smallest
non-zero eigenvalue of the Hodge Laplacian on $p$-forms for a
compact Riemannian manifold with bounded sectional curvature,
while the purpose of \cite{MG} is to study the spectrum of
$\Delta$ on compact hyperbolic 3-dimensional manifolds. In
particular, McGowan develops in \cite{MG} a quite general method
to bound from below "small" eigenvalues of $\Delta$ on compact
manifolds (Lemma 2.3 in \cite{MG}).

\smallskip

The purpose of this paper is in some sense to improve or to unify
these results in the context given by the discretization. More
precisely, if $M$ is a compact Riemannian manifold and if $X$ is a
discretization of $M$ (in the sense of \cite{Ch}), we obtain
naturally from $X$ a finite open cover $\U_X$ which will be
contractible if the mesh of the discretization is sufficiently
small. To such an open cover we can associate the complex of \Cech
cochains naturally endowed with a coboundary operator $\delta$.
Moreover, with an inner product on \Cech cochains, we can
construct the adjoint of $\delta$, namely $\delta^*$ and define
the following combinatorial Laplacian $\check{\Delta} = \delta
\delta^* + \delta^* \delta$.

\smallskip

The main result consists in establishing a uniform comparison
between the spectrum of the Hodge Laplacian and the spectrum of
such a combinatorial Laplacian. That is to say, if $\K$ denotes
the set of compact connected Riemannian manifolds with bounded (by
$\kappa$) sectional curvature and injectivity radius bounded from
below by $r_0$, we show that there exists a positive constant
$\rho_0$ depending only on $n$, $\kappa$ and $r_0$ such that if we
fix $0 < 3\eps < \rho_0$, there exist positive constants $c_1$ and
$c_2$ depending only on $n$, $p$, $\kappa$ and $\eps$ such that
for any $M \in \K$ and for any $\eps$-discretization $X$ of $M$ we
can compare the $k^{\text{th}}$ eigenvalue of $\Delta$ on
$p$-forms to the $k^{\text{th}}$ eigenvalue of $\check{\Delta}$ on
\Cech $p$-cochains (for $1 \leq p \leq n-1$) in the following way
\begin{eqnarray*}
c_1 \lambda_{k,p}(X) \leq \lambda_{k,p}(M) \leq c_2
\lambda_{k,p}(X)
\end{eqnarray*}
for any $k \leq K$ and $K$ depends on $X$ (see Theorem
\ref{thmmain} for the precise statement).

\smallskip

As an application of Theorem \ref{thmmain}, we obtain a lower
bound for the first non-zero eigenvalue of $\Delta$ (see Theorem
\ref{lambda1}) in terms of the volume of the manifold. This result
has to be compared with the result obtained by Chanillo and Trèves
(Theorem 1.1, in \cite{CT}). In their proof, the authors use in a
crucial manner a lemma due to Trèves (Lemma A.5 in \cite{Tr})
which turns out to be false (see Remark \ref{faute}). In Lemma
\ref{lmmT}, we state and prove a "weaker" version of Trèves'
lemma. A direct corollary of this lemma is a lower bound for the
spectrum of the combinatorial Laplacian (see Theorem
\ref{lambda2}) and so, thanks to Theorem \ref{thmmain}, a lower
bound for the spectrum of $\Delta$ (see Theorem \ref{lambda1}).

\smallskip

As another consequence of the proof of Theorem \ref{thmmain}, we
obtain a version of McGowan's lemma (Lemma 2.3 in \cite{MG})
slightly more general as it is concerned with $p$-forms on compact
Riemannian manifolds with bounded sectional curvature, but not so
general as it is valid only for contractible open covers (see
Lemma \ref{lmmmcg}). Finally, another interesting application of
the method developed here concerns Whitney forms. Indeed, Whitney
forms come out in \cite{DP} as a natural way to smooth \Cech
cochains. Nevertheless, in order to keep a uniform comparison of
the spectra, the results given in \cite{DP} on Whitney forms are
not useful to our purpose. Hence, we obtain as a corollary of the
method, the appropriate results to show that Whitney forms are
even so a suitable tool to smooth \Cech cochains (see Section
\ref{secW}).

\smallskip

The paper is organized as follows. In Section \ref{sec2}, we begin
by recalling different definitions and properties of differential
forms and \Cech cochains. In particular, in Section \ref{pfdeR},
we sketch the proof of the de Rham Theorem due to A. Weil as it
will be the starting point of the proof of Theorem \ref{thmmain}.
Finally, we recall the definition of a discretization and its main
properties.

\smallskip

Section \ref{sec3} is devoted to the proof of Theorem
\ref{thmmain}. The basic idea of the proof is to associate a \Cech
cochain to a differential form via a discretizing operator and
vice versa via a smoothing operator, in order to compare "small"
eigenvalues. These operators are essentially constructed as in the
proof (of A. Weil) of the de Rham Theorem thanks to the \Cech - de
Rham double complexe. To that aim, we need a few technical
results. In particular, we need a normed version of the Poincaré
Lemma and a similar result for \Cech cochains. This is done in
Lemma \ref{poincare2} and in Lemma \ref{poinc2}. Moreover, as in
\cite{MG}, it is necessary to bound from below the spectrum of
$\Delta$ with absolute boundary conditions on finite intersections
of open sets of the open cover. To that aim, we show that for a
sufficiently small $\eps$, the intersection of balls of radius
$\eps$ is convex and is quasi-isometric to a Euclidean convex.
Thanks to a result of Guerini (\cite{Gu}) we can then bound from
below the spectrum of such intersections (this appears in Section
\ref{sec2} as properties of the discretization, see Lemma
\ref{lmmrho} and Lemma \ref{lmmmu}). Note that Chanillo and Trèves
met also this problem and they solve it using a (finite) sequence
of open covers and with Lemma 2.2 in \cite{CT} (which is a
consequence of a normed version of the Poincaré Lemma in the
Euclidean setting). For "large" eigenvalues, it suffices to have
an upper bound for the $k^{\text{th}}$ eigenvalue of $\Delta$ and
of $\check{\Delta}$ to have the claim.

\smallskip

In Section \ref{sec4}, we present the consequences of Theorem
\ref{thmmain} mentioned above.

\smallskip

Finally, in the appendix we recall the (more or less classical)
definition and the properties of Whitney forms. At the end of the
appendix, we give the proof of the technical lemma about the
Euclidean convexity of the intersection of small balls.

%% file: DiscFormsSettings.tex
\section{Settings} \label{sec2}

In this section, we recall some definitions and basic facts on the
Laplacian acting on differential forms and on the Laplacian acting
on \Cech cochains. For the convenience of the reader and as it is
a key tool for the paper, a paragraph is also devoted to the
sketch of a classical proof due to A. Weil of the de Rham Theorem
(for contractible open covers) relying on the \Cech - de Rham
double complexe (see for instance Appendix A of \cite{Go} or
Chapter 3 of \cite{Mo}). Finally, we define the discretization of
a manifold and discuss some of its properties.


\subsection{Laplacian acting on differential forms} \label{laplF}

Let $(M^n,g)$ be a compact connected $n$-dimensional Riemannian
manifold without boundary. Denote by $\Lambda^p(M)$ the vector
space of smooth differential $p$-forms, for $0 \leq p \leq n$. Let
$d : \Lambda^p(M) \to \Lambda^{p+1}(M)$ be the exterior
differential and $d^* : \Lambda^{p+1}(M) \to \Lambda^p(M)$ its
formal adjoint (with respect to the $L^2$-inner product) the
codifferential. Then the Laplacian acting on $p$-forms is defined
by $\Delta : \Lambda^{p}(M) \to \Lambda^p(M)$, $\Delta = d d^* +
d^*d$. The spectrum of $\Delta$ is discrete and will be denoted by
$$0 < \lambda_{1,p}(M) \leq \lambda_{2,p}(M) \leq \ldots \leq
\lambda_{k,p}(M) \leq \ldots$$ where $0$ is of multiplicity
$b_p(M)$ and the positive eigenvalues are repeated as many times
as their multiplicity. Let us recall that half of the spectrum is
redundant. That is to say, if $\lambda > 0$ is an eigenvalue of
$\Delta$ on $p$-forms and if $E_p(\lambda)$ denotes the
$\lambda$-eigenspace, then $E_p(\lambda)$ splits as follows
$E_p(\lambda) = E_p^{d^*}(\lambda)\oplus E_p^{d}(\lambda)$ where
$E_p^{d^*}(\lambda)=\{\omega \in E_p(\lambda) : d^*\omega = 0\}
\subseteq d^*\Lambda^{p+1}(M)$ is the $\lambda$-eigenspace of
$d^*d$ and $E_p^{d}(\lambda)=\{\omega \in E_p(\lambda) : d\omega =
0\} \subseteq d\Lambda^{p-1}(M)$ is the $\lambda$-eigenspace of $d
d^*$. Moreover, $d^*$ maps $E_p^d(\lambda)$ isomorphically onto
$E_{p-1}^{d^*}(\lambda)$ and $d$ maps $E_p^{d^*}(\lambda)$
isomorphically onto $E_{p+1}^{d}(\lambda)$. Hence, $E_p(\lambda) =
E_p^{d^*}(\lambda) \oplus E_{p-1}^{d^*}(\lambda)$. So for our
purpose it will be sufficient to study the spectrum of $d^*d$ on
coexact forms.

\smallskip

Let $\lambda_{k,p}^{d^*}(M)$ the $k^{\text{th}}$ (positive)
eigenvalue of $d^*d : d^*\Lambda^{p+1}(M) \to
d^*\Lambda^{p+1}(M)$. The following variational characterization
of the spectrum of $d^*d$ holds $$\lambda_{k,p}^{d^*}(M) =
\min\limits_{\Sigma^k}\max\left\{
\frac{\|d\omega\|^2}{\|\omega\|^2} : \omega \in \Sigma^k \setminus
\{0\} \right\}$$ where $\Sigma^k$ ranges over all $k$-dimensional
vector subspaces of $d^*\Lambda^{p+1}(M)$ and $\|\cdot\|$ denotes
the $L^2$-norm for differential forms.


\subsection{\Cech cohomology and combinatorial
Laplacian}\label{laplC}

Let $M^n$ be a compact connected $n$-dimensional manifold. Let $\U
= \{U_i\}_{1\leq i \leq N}$ be a finite open cover of $M$. The
nerve of $\U$, denoted by $N(\U)$, is the simplicial complex whose
set of $q$-simplices is given by $$S_q(\U) = \{(i_0,\ldots,i_q) :
i_0 < \ldots < i_q\text{ and } U_{i_0}\cap \ldots \cap U_{i_q}
\neq \emptyset\}$$ for any $q \geq 0$. A \Cech $q$-cochain is an
application $c : S_q(\U) \to \R$. Denote by $\CqR{q}{\U}$ the set
of \Cech $q$-cochains. Let us remark that $\CqR{q}{\U}$ is
naturally endowed with a vector space structure and let us define
a coboundary operator $\delta : \CqR{q}{\U} \to \CqR{q+1}{\U}$ by
$$\delta c (i_0,\ldots,i_{q+1}) = \sum\limits_{j=0}^{q+1}
(-1)^jc(i_0,\ldots,i_{j-1},i_{j+1},\ldots,i_{q+1})$$ for any
$\{i_0,\ldots,i_{q+1}\} \in S_{q+1}(\U)$. Then $\delta \circ
\delta = 0$ and the cochain complex $\{\CqR{q}{\U}, \delta\}$
gives rise to the \Cech cohomology groups of the cover $\U$,
$\check{H}^*(\U)$.

Endow then $\CqR{q}{\U}$ with the following scalar product, for
any $c_1$, $c_2 \in \CqR{q}{\U}$ $$(c_1,c_2) = \sum_{I \in
S_q(\U)}c_1(I) c_2(I)$$ and consider $\delta^* : \CqR{q+1}{\U} \to
\CqR{q}{\U}$ the adjoint of $\delta$ with respect to
$(\cdot,\cdot)$.

\begin{dfn} The \textbf{combinatorial Laplacian} $\check{\Delta} :
\CqR{q}{\U}  \to \CqR{q}{\U}$ is defined by $\check{\Delta} =
\delta \delta^* + \delta^* \delta$.
\end{dfn}

The combinatorial Laplacian is self-adjoint and non-negative by
definition. Its spectrum will be denoted by $$0 <
\lambda_{1,q}(\U) \leq \lambda_{2,q}(\U)\leq \ldots \leq
\lambda_{L,q}(\U)$$ where $0$ is of multiplicity $\check{b}_q(\U)$
and $L + \check{b}_q(\U) = dim (\CqR{q}{\U}) = |S_q(\U)|$. As for
the Laplacian on differential forms, half of the spectrum is
redundant i.e. if $\lambda > 0$ is an eigenvalue of
$\check{\Delta}$ on \Cech $q$-cochains and if
$\check{E}_q(\lambda)$ denotes the $\lambda$-eigenspace, then
$\check{E}_q(\lambda) = \check{E}_q^{\delta^*}(\lambda) \oplus
\check{E}_{q-1}^{\delta^*}(\lambda)$ where
$\check{E}_q^{\delta^*}(\lambda)$ is the $\lambda$-eigenspace of
$\delta^* \delta$ acting on $\delta^* \CqR{q+1}{\U}$. So for our
purpose it will be sufficient to study the spectrum of $\delta^*
\delta$ on $\delta^* \CqR{q+1}{\U}$ i.e. on coexact \Cech
cochains. In the sequel, $\lambda^{\delta^*}_{k,q}(\U)$ denotes
the $k^{\text{th}}$ (positive) eigenvalue of $\delta^* \delta :
\delta^{*} \CqR{q+1}{\U} \to \delta^*\CqR{q+1}{\U}$. The following
variational characterization holds $$\lambda_{k,q}^{\delta^*}(\U)
= \min_{V^k} \max\left\{\frac{\|\delta c\|^2}{\|c\|^2} : c \in V^k
\setminus\{0\} \right\}$$ where $V^k$ ranges over all
$k$-dimensional vector subspaces of $\delta^*\CqR{q+1}{\U}$.


\subsection{De Rham Theorem}\label{pfdeR}

Recall that an open cover $\U$ is called contractible if for any
$I \in S_q(\U)$, $U_I = \bigcap_{i \in I}U_i$ is contractible. The
following theorem is due to de Rham.

\begin{thm}
Let $(M^n,g)$ be a compact connected $n$-dimensional Riemannian
manifold without boundary. Let $\U$ be a contractible finite open
cover of $M$. Then the $p^{\text{th}}$ group of de Rham's
cohomology  $H^p(M)$ is isomorphic to $\check{H}^p(\U)$.
\end{thm}

\begin{rmk} Note that a consequence of the de Rham Theorem is that if $\U$ is a
contractible cover, then $b_p(M) = \check{b}_p(\U)$.
\end{rmk}

Let us introduce now the vector spaces $\Cq{q}{U}{\Lambda^p}$ of
$q$-cochains of $p$-forms i.e. $c$ is in $\Cq{q}{U}{\Lambda^p}$ if
$c(I)$ is a $p$-form on $U_I$ for any $I$ in $S_q(\U)$ . Define
then the following coboundary operators
\begin{multline}\nonumber \delta : \Cq{q}{U}{\Lambda^p} \to
\Cq{q+1}{U}{\Lambda^p}\text{ defined by } \\ \delta
c(i_0,\ldots,i_{q+1}) = \sum_{j=0}^{q+1}(-1)^j
c(i_0,\ldots,i_{j-1},i_{j+1},\ldots,i_{q+1})
\end{multline} for any
$\{i_0,\ldots,i_{q+1}\}\in S_{q+1}(\U)$ and $$d :
\Cq{q}{U}{\Lambda^p} \to \Cq{q}{U}{\Lambda^{p+1}} \text{ defined
by } dc(I) = d(c(I))$$ for any $I \in S_q(\U)$. Then $d \circ d
=0$, $\delta \circ \delta = 0$ and $d \circ \delta = \delta \circ
d$.  The \Cech- de Rham double complex is the following
commutative diagram,
\begin{figure}
\begin{equation} \nonumber
   \scalebox{.68}{
   \begin{CD}   & & \CqR{0}{\U} @> \delta >> \CqR{1}{\U} @>\delta>> \ldots
    @>\delta >> \CqR{q-1}{\U} @>\delta >>  \CqR{q}{\U} @>\delta >>\ldots   \\
     & & @VViV                      @VViV                   & &    @VViV                       @VViV                    \\
     \Lambda^{0}(M)  @> r >> \Cq{0}{U}{\Lambda^0} @> \delta >>      \Cq{1}{U}{\Lambda^0} @> \delta >> \ldots
@>\delta >> \Cq{q-1}{U}{\Lambda^0} @>\delta >>
\Cq{q}{U}{\Lambda^0} @>\delta >>\ldots   \\
     @VVdV @VVdV                      @VVdV                   & & @VVdV                          @VVdV                     \\
     \Lambda^{1}(M)  @> r >> \Cq{0}{U}{\Lambda^1} @> \delta >>      \Cq{1}{U}{\Lambda^1} @> \delta >> \ldots
@>\delta >> \Cq{q-1}{U}{\Lambda^1} @>\delta >>
\Cq{q}{U}{\Lambda^1} @>\delta >>\ldots   \\
     @VVdV @VVdV                      @VVdV                   & & @VVdV                          @VVdV                    \\
        \vdots && \vdots && \vdots && \cdots && \vdots && \vdots\\
     @VVdV @VVdV                      @VVdV                   & &      @VVdV                     @VVdV                   \\
     \Lambda^{p-1}(M)  @> r >> \Cq{0}{U}{\Lambda^{p-1}} @> \delta >>      \Cq{1}{U}{\Lambda^{p-1}} @> \delta >> \ldots
     @>\delta >> \Cq{q-1}{U}{\Lambda^{p-1}} @>\delta >>  \Cq{q}{U}{\Lambda^{p-1}} @>\delta >> \ldots   \\
     @VVdV @VVdV                      @VVdV                   & &           @VVdV                @VVdV                    \\
     \Lambda^{p}(M)  @> r >> \Cq{0}{U}{\Lambda^p} @> \delta >>      \Cq{1}{U}{\Lambda^p} @> \delta >> \ldots
     @>\delta >> \Cq{q-1}{U}{\Lambda^p} @>\delta >> \Cq{q}{U}{\Lambda^p} @>\delta >>\ldots   \\
     @VVdV @VVdV                      @VVdV                   & &          @VVdV                 @VVdV                    \\
        \vdots && \vdots && \vdots &&  && \vdots & & \vdots
    \end{CD}}
\end{equation}

\caption{The \Cech - de Rham double complexe.}
\end{figure}
where $r$ denotes the restriction map to each open of the cover
and $i$ the natural injection. The first step in the proof of the
de Rham Theorem is to show that the rows (except the first) and
the columns (except the first) of this diagram are exact. This is
a direct consequence of the Poincaré Lemma (Lemma \ref{poincare})
and Lemma \ref{poinc}.

\begin{lmm}\label{poincare}
Let $p > 0$. Let $\U$ be a contractible cover. Let $\omega \in
\Cq{q}{U}{\Lambda^p}$ such that $d\omega = 0$. Then there exists
$\eta \in \Cq{q}{U}{\Lambda^{p-1}}$ such that $d \eta = \omega$.
\end{lmm}
\textbf{Proof}: see \cite{Go}, A.6. $\square$

\begin{lmm}\label{poinc}
Let $q > 0$. Let $c \in \Cq{q}{U}{\Lambda^p}$ such that $\delta c
= 0$. Then there exists $b \in \Cq{q-1}{U}{\Lambda^p}$ such that
$\delta b = c$.
\end{lmm}
\textbf{Proof}: see \cite{Go}, proof of Lemma A.4.1. $\square$

\smallskip

The proof of the de Rham Theorem goes then as follows. Let $\omega
\in \Lambda^p(M)$ such that $d \omega =0$. Let $f_0 = r(\omega)
\in \Cq{0}{\U}{\Lambda^p}$, then $df_0 = 0 = \delta f_0$ and the
system of equations $$f_0 = df_1\ ,\ \delta f_1 = d f_2\ , \
\delta f_2 = d f_3\ ,\ \ldots \ , \ \delta f_{p-1} = d f_p$$ has a
solution with $f_j \in \Cq{j-1}{U}{\Lambda^{p-j}}$ for $j \geq 1$.
Moreover, $\delta (\delta f_p) = 0$, hence $\delta f_p \in
\CqR{p}{\U}$. The application $\Psi : \{\omega \in \Lambda^p(M) :
d\omega = 0 \} \to \{ c \in \CqR{p}{\U} : \delta c =0 \}$ given by
$\Psi(\omega) = \delta f_p$, where $f_p$ is constructed as above,
induces an isomorphism in cohomology. In particular, if $\omega$
is exact, $\Psi(\omega)$ is also exact i.e. there exists $c \in
\CqR{p-1}{\U}$ such that $\delta c = \Psi(\omega)$ (note that in
general $f_p \notin  \CqR{p-1}{\U}$). Naturally, we can construct
another application going from closed \Cech $p$-cochains to closed
$p$-forms exactly in the same way and obtain also an isomorphism
in cohomology. $\square$



\subsection{Discretization of a manifold}

Let $(M^n,g)$ be a connected compact $n$-dimensional Riemannian
manifold without boundary. Let $\eps > 0$.

\begin{dfn} An $\eps$-discretization $X$ of $M$ is a maximal
$\eps$-separated subset of $M$ i.e. $X$ is a subset of $M$
satisfying
\begin{enumerate}
\item[(i)] $\forall p \neq q \in X$, $d(p,q) \geq \eps$,
\item[(ii)] $\U_X = \{ B(p,\eps)\}_{p\in X}$ is an open cover
of $M$.
\end{enumerate}
\end{dfn}
Note that as $M$ is compact, $X$ is finite of cardinality $|X|$.
So we can number the elements of $X = \{p_1,\ldots,p_{|X|}\}$ and
denote $U_i = B(p_i,\eps)$, for $i=1,\ldots,|X|$. In particular,
any discretization of $M$ gives rise to a combinatorial Laplacian
$\check{\Delta}$ as defined in Section \ref{laplC}. In the sequel,
$\lambda_{k,q}(X)$ will denote the $k^{\text{th}}$ eigenvalue of
the combinatorial Laplacian associated to the open cover $\U_X$
acting on \Cech $q$-cochains i.e. $\lambda_{k,q}(X) =
\lambda_{k,q}(\U_X)$.

\smallskip

Note also that if $\eps$ (the mesh of the discretization) is
smaller than the convexity radius of $M$, then $\U_X$ is a
contractible open cover and $\check{b}_p(\U_X) = b_p(M)$.

\begin{dfn} For $\kappa \geq 0$, $r_0 > 0$ and $n \in \N^*$, we
define $\K$ as the set of all connected compact $n$-dimensional
Riemannian manifold $(M^n,g)$ without boundary with uniformly
bounded sectional curvature i.e. $|K_g| \leq \kappa$ and
injectivity radius bounded below i.e. $Inj(M,g) \geq r_0$.
\end{dfn}

\begin{rmk}\label{cardX} For $n\in \N^*$, $\kappa\geq 0$, $r_0 > 0$ and $0 < 2\eps < r_0$,
there exists $\nu(n,\kappa) >0$ such that, for any $(M,g) \in \K$
and any $\eps$-discretization $X$ of $M$, the cardinality of $\{j
: U_j \cap U_I \neq \emptyset \}$ is bounded above by $\nu$, for
any $I \in S_q(\U_X)$. This is a direct consequence of the
Bishop-Gromov volume comparison Theorem (see for instance
\cite{Ch}, Lemma V.3.1, p.147). Furthermore, by Croke's Inequality
and Bishop's comparison Theorem (see \cite{Ch} p.126 and p.136) we
can assert that there exist positive constants $c_1$, $c_2$
depending only on $n$, $\kappa$ and $\eps$ such that $c_1 Vol(M)
\leq |X| \leq c_2 Vol(M)$. In particular, we obtain that
$|\S_q(\U_X)| \leq \frac{\nu^{q}}{(q+1)!} |X| \leq
\frac{\nu^{q}}{(q+1)!} c_2 Vol(M)$.
\end{rmk}

The following lemma shows that in general a sufficiently small
ball is quasi-isometric (in the sense of \cite{Do}, (3.2)) to a
Euclidean convex. In particular, this will imply that on
intersections of sufficiently small balls we can find a lower
bound for the first positive eigenvalue of $\Delta$ with absolute
boundary condition (see Lemma \ref{lmmmu}). This is an essential
result for the discretization as we will see later.

\begin{lmm}\label{lmmrho}
Let $n \in \N^*$, $\kappa \geq 0$ and $r_0 > 0$. There exists a
constant $0 < \rho_0 < r_0$ depending only on $n$, $\kappa$ and
$r_0$ such that for any $(M,g) \in \K$ and for any $p \in M$,
there exist a Euclidean convex $C_p \subseteq \R^n$ and a
diffeomorphism $\varphi : C_p \to B(p,\rho_0)$ such that for any
$B(q,\rho) \subseteq B(p,\rho_0)$, the ball $B(q,\rho)$ is convex
and $\varphi^{-1}(B(q,\rho))$ is a Euclidean convex. Moreover,
$(B(q,\rho),g)$ is quasi-isometric to $B(q,\rho)$ endowed with the
Euclidean metric induced by $\varphi^{-1}$ and the constants of
quasi-isometry depend only on $n$, $\kappa$ and $d(p,q) + \rho$.
\end{lmm}
\textbf{Proof}: see Appendix \ref{pflmmrho}. $\square$

\smallskip

Note that the intersection of small balls is a convex with not
necessarily smooth boundary. So that it is not obvious that in
this case the spectrum of the Laplacian with absolute boundary
condition is discrete. In \cite{MiMi}, the authors show that the
spectrum of the Laplacian with absolute (or relative) boundary
condition is discrete even if the boundary is only given by a
Lipschitz function (Proposition 5.3 in \cite{MiMi}). Moreover,
Theorem 5.1 of \cite{MTV} implies that the following classical
variational characterization of the spectrum is still valid for
bounded convex domains i.e. if $\Omega$ is a bounded convex domain
of $M$, then the $k^{\text{th}}$ eigenvalue of the Laplacian for
$p$-forms on $\Omega$ with absolute boundary condition is given by
$$\lambda_{k,p}^{abs}(\Omega) = \min\limits_{\Sigma^k}\max\left\{
\frac{\|d\omega\|^2 + \|\delta \omega\|^2}{\|\omega\|^2} : \omega
\in \Sigma^k \setminus \{0\} \text{ such that } i_\nu(\omega) = 0
\right\}$$ where $\Sigma^k$ ranges over all $k$-dimensional vector
subspaces of $\Lambda^{p}(\Omega)$ and $i_\nu$ is the interior
product by $\nu$ the outward pointing normal unit vector to the
boundary (defined almost everywhere). In particular, the result on
quasi-isometric metrics of Dodziuk (Proposition 3.3 of \cite{Do})
is valid in this context.

\begin{lmm}\label{lmmmu}
Let  $n \geq 2$, $\kappa \geq 0$, $r_0 > 0$ and let $\rho_0$ given
by Lemma \ref{lmmrho}. Let $0 < 3 \eps < \rho_0$. Then there
exists a positive constant $\mu(n,\kappa,\eps)$ depending only on
$n$, $\kappa$ and $\eps$ such that for any $(M,g) \in \K$ and for
any $\eps$-discretization $X$ of $M$ $$\lambda_{1,p}^{abs}(U_I)
\geq \mu(n,\kappa,\eps)$$ for any $p = 0,\ldots,n$ and any $I \in
S_q(\U_X)$, $q\geq 0$.
\end{lmm}
\textbf{Proof}: let $(M,g) \in \K$ and $X$ an
$\eps$-discretization of $M$ with $0< 3 \eps < \rho_0$. Fix $p \in
X$ and let $q \in X$ such that $B(p,\eps) \cap B(q,\eps) \neq
\emptyset$. Then $B(q,\eps) \subseteq B(p,3 \eps) \subseteq
B(p,\rho_0)$. By Lemma \ref{lmmrho}, there exists a diffeomorphism
$\varphi$ such that $\varphi^{-1}(B(q,\eps))$ is a Euclidean
convex for any $q \in X$ such that $B(q,\eps) \cap B(p,\eps) \neq
\emptyset$. In particular, $\varphi^{-1}\left( B(p,\eps)\cap
B(q,\eps)\right)$ is an intersection of Euclidean convexes and as
such it is a Euclidean convex. Moreover, $\varphi^{-1}$ restricted
to $B(p,3\eps)$ is a quasi-isometry with constants of
quasi-isometry depending only on $n$, $\kappa$ and $\eps$. Let
$U_I$ a non-empty finite intersection of elements of $\U_X$ and
$V_I = \varphi^{-1}(U_I)$ the Euclidean convex which is
quasi-isometric to $U_I$ via $\varphi$ i.e.
$(\varphi(V_I),(\varphi^{-1})^*(eucl))$ is quasi-isometric to
$(U_I,g)$ with constants of quasi-isometry $\alpha$ depending only
on $n$, $\kappa$ and $\eps$ (i.e. $\alpha ^{-1}
(\varphi^{-1})^*(eucl) \leq g \leq \alpha(\varphi^{-1})^*(eucl)$).
Then by Proposition 3.3 of \cite{Do}, there exist positive
constants $c_1$ and $c_2$ depending only on $\alpha$ and $n$ such
that
\begin{eqnarray}\label{rido}
c_1 \lambda_{1,p}^{abs}(U_I,(\varphi^{-1})^*(eucl)) \leq
\lambda_{1,p}^{abs}(U_I,g) \leq
c_2\lambda_{1,p}^{abs}(U_I,(\varphi^{-1})^*(eucl)) .
\end{eqnarray}
Note that $(U_I,(\varphi^{-1})^*(eucl))$ is a Euclidean convex of
diameter bounded above by $d(n,\kappa,\eps)$. Finally, Guerini
shows in \cite{Gu}, that the first eigenvalue of the Laplacian
with absolute boundary condition on a Euclidean convex with smooth
boundary is bounded below by a constant depending on the diameter
of the convex. Note that Guerini's proof can be adapted
straightforward to obtain the same result for convexes with
piecewise smooth boundary. Hence, we obtain that there exists a
positive constant $c(n,p)$ such that
\begin{eqnarray}\label{lgu}
\lambda_{1,p}^{abs}(U_I,(\varphi^{-1})^*(eucl)) \geq
\frac{c(n,p)}{diam(U_I,(\varphi^{-1})^*(eucl))^2} \geq
\frac{c(n,p)}{d(n,\kappa,\eps)^2}
\end{eqnarray}
Finally, (\ref{rido}) and (\ref{lgu}) imply the claim. $\square$


%% file: DiscFormsMain.tex
\section{Comparison of spectra} \label{sec3}

This section is devoted to the proof of the main theorem of the
paper. Let us state the result.

\begin{thm}\label{thmmain}
Let $n \geq 2$, $\kappa \geq 0$, $r_0 > 0$. Let
$\rho_0(n,\kappa,r_0)$ be given by Lemma \ref{lmmrho} and $0 < 3
\eps < \rho_0$. Let $1 \leq p \leq n-1$. Then there exist positive
constants $c_1$, $c_2$ depending only on $n$, $p$, $\kappa$ and
$\eps$ such that for any $M \in \K$ and for any
$\eps$-discretization $X$ of $M$, we have $$c_1 \lambda_{k,p}(X)
\leq \lambda_{k,p}(M) \leq c_2 \lambda_{k,p}(X)$$ for any $1 \leq
k \leq |\CqR{p}{\U_X}| - \check{b}_p(\U_X) = |\CqR{p}{\U_X}| -
b_p(M)$.
\end{thm}

As we have seen before (in Section \ref{laplF}), it will be
sufficient to establish the result for the spectrum of $d^*d$ on
coexact $p$-forms and for the spectrum of $\delta^*\delta$ on
coexact \Cech $p$-cochains. The proof goes in two steps. First
step consists in comparing "small" eigenvalues. We need to
construct a discretizing operator that associates to a coexact
$p$-form a coexact \Cech $p$-cochain (see Section \ref{constrD})
and a smoothing operator that goes in the opposite direction (see
Section \ref{constrS}), in order to compare their respective
Rayleigh quotients. The idea is to proceed as in the proof of the
de Rham Theorem and use the \Cech - de Rham double complexe. But
as we need a control of the norms involved, we have to establish
versions of the Poincaré Lemma (Lemma \ref{poincare}) and of Lemma
\ref{poinc} with a suitable control of the norms (see Lemma
\ref{poincare2} and Lemma \ref{poinc2}). The second step of the
proof deals with "large" eigenvalues and is reduced to find upper
bounds for the $k^{\text{th}}$ eigenvalues involved depending only
on the parameters of the problem (see Section \ref{secup}).

\smallskip

In the sequel, we consider $(M,g)$ in $\K$ and $X$ an
$\eps$-discretization with $0 < 3 \eps < \rho_0$. Denote by $\U$
the open cover induced by $X$ i.e. $\U = \{U_i = B(p_i,\eps) :
i=1,\ldots,|X|\}$ and fix $1 \leq p \leq n-1 $.

\subsection{From smooth forms to \Cech cochains}\label{constrD}

In this section, we are going to construct $$\D : d^*
\Lambda^{p+1}(M) \to \delta^*\CqR{p+1}{\U}$$ such that there exist
positive constants $c_1$, $c_2$ and $\Lambda$ depending only on
$n$, $p$, $\kappa$ and $\eps$ such that
\begin{enumerate}
    \item[$(i)_{\D}$] $\|\delta \D (\omega) \|^2 \leq c_1 \|d
    \omega \|^2$, for any $\omega \in d^* \Lambda^{p+1}(M)$,
    \item[$(ii)_{\D}$] $\|\D \omega \|^2 \geq c_2 \|\omega\|^2$,
    for any $\omega \in d^* \Lambda^{p+1}(M)$ satisfying $\|
    d\omega \|^2 \leq \Lambda \| \omega\|^2$.
\end{enumerate}
To that aim, we need the following version of the Poincaré Lemma.
Note that this lemma will be verified in particular by any
non-empty intersection of open sets in $\U$ thanks to Lemma
\ref{lmmmu} (where $\mu$ depends on $n$, $\kappa$ and $\eps$).
\begin{lmm}\label{poincare2}
    Let $U$ be a contractible open set such that $\lambda_{1,p}^{abs,d}(U)
    \geq \mu >0$, ($1\leq p\leq n$). Let $\omega$ be a closed
    $L^2$-integrable $p$-form on $U$ i.e. $d\omega = 0$. Then
    there exists $\eta \in \Lambda^{p-1}(U)$ such that $d\eta =
    \omega$ and $\|\eta\|^2_{L^2(U)} \leq \frac{2}{\mu}
    \|\omega\|^2_{L^2(U)}$.
\end{lmm}
\textbf{Proof}: we have the following characterization of the
first eigenvalue of the Laplacian on exact $p$-forms (see
Proposition 3.1. of \cite{Do} or Proposition 2.1. of \cite{MG}),
$$\lambda_{1,p}^{abs,d}(U) = \inf_{V} \sup \left\{
\frac{\|\omega\|^2_{L^2(U)}}{\|\eta\|^2_{L^2(U)}} : \omega \in
V\setminus\{0\} \text{ , } d\eta = \omega  \right\}$$ where $V$
ranges over all $1$-dimensional vector subspaces of exact
$p$-forms. If $\omega \in \Lambda^p(U)$ is closed, by the Poincaré
Lemma $\omega$ is exact. So that we get $$ \mu \leq
\lambda_{1,p}^{abs,d}(U) \leq \sup \left\{
\frac{\|\omega\|^2_{L^2(U)}}{\|\eta\|^2_{L^2(U)}} :  d\eta =
\omega  \right\}$$ and hence there exists $\eta \in
\Lambda^{p-1}(U)$ such that $d\eta = \omega$ and $\frac{1}{2} \mu
\leq \frac{\|\omega\|^2_{L^2(U)}}{\|\eta\|^2_{L^2(U)}}$ which is
the claim. $\square$

\smallskip

\begin{rmk}\label{norm}
Let us introduce the following norm. If $c \in
\Cq{q}{U}{\Lambda^p}$ let $$\|c\|^2 = \sum_{I \in S_q(\U)}
\|c(I)\|^2_{L^2(U_I)}$$ where $\|\cdot\|_{L^2(U_I)}$ denotes the
$L^2$-norm for $p$-forms on $U_I$. In particular, if $\omega$ is a
$p$-form on $M$ and $r$ is the restriction to each open of $\U$,
then there exist positive constants $c_1$ and $c_2$ depending only
on $n$, $\kappa$ and $\eps$ such that $c_1 \|r(\omega)\|^2 \leq
\|\omega\|^2 \leq c_2 \|r(\omega)\|^2$.
\end{rmk}

\subsubsection*{Construction by induction of $\D$}

Let $\omega \in d^*\Lambda^{p+1}(M)$. The goal is to construct $\D
(\omega) \in \delta^*\CqR{p+1}{\U}$. The idea is to consider $d
\omega$ which is an exact $(p+1)$-form and to construct an exact
\Cech $(p+1)$-cochain $\delta \D(\omega)$ such that $(i)_{\D}$
holds. A suitable candidate for $\delta \D (\omega)$ is the \Cech
cochain given by the proof of the de Rham Theorem and the double
complexe. Moreover, the double complexe and the normed version of
the Poincaré Lemma give almost directly the inequality $(i)_{\D}$,
whereas $(ii)_{\D}$ is not a so direct consequence of the
construction. Hence, as suggested in \cite{CT}, we construct an
auxiliary $p$-form thanks to Whitney forms to obtain $(ii)_{\D}$.
We proceed by induction.

\smallskip

\textbf{First step of induction}: define $c_{p+1,0} \in
\Cq{0}{U}{\Lambda^{p+1}}$ by $c_{p+1,0} = r(d\omega)$ i.e.
$c_{p+1,0}(i) = d\omega_{|_{U_i}}$. Then $dc_{p+1,0} = 0 = \delta
c_{p+1,0}$ and $W(c_{p+1,0}) = d\omega$, where $W$ is the Whitney
map defined in Appendix \ref{whitney}. Then there exist positive
constants $c_1$, $c_2$ and $c_3$ depending only on $n$, $p$,
$\kappa$ and $\eps$ such that the three following assertions hold.
\begin{enumerate}
\item[$(a)_1$] There exists $c_{p,0} \in \Cq{0}{U}{\Lambda^{p}}$
such that $d c_{p,0} = c_{p+1,0}$ and $\|c_{p,0}\|^2 \leq c_1 \| d
\omega \|^2$.
\item[$(b)_1$] Let $c_{p,1} = \delta c_{p,0}$. We have $\delta
c_{p,1} = 0 = dc_{p,1}$ and $\|c_{p,1}\|^2 \leq c_2
\|d\omega\|^2$.
\item[$(c)_1$] Let $v^{(1)} = W(c_{p,0}) \in \Lambda^{p}(M)$. We
have $dv^{(1)} = d\omega + W(c_{p,1})$ and $\|v^{(1)}\|^2 \leq c_3
\|d\omega\|^2$.
\end{enumerate}
Indeed, $(a)_1$ is a direct consequence of Lemma \ref{poincare2},
of the definition of $c_{p+1,0}$ and of Remark \ref{norm}. Then,
clearly $\delta c_{p,1} = 0$ and $d c_{p,1} = \delta d c_{p,0} =
\delta c_{p+1,0} = 0$.  Moreover, there exists $c(n,\kappa,\eps)$
such that for any cochain $\|\delta b\|^2 \leq c \|b\|^2$ (see
(\ref{updel})) and combined with $(a)_1$ this implies $(b)_1$.
Finally, by Lemma \ref{W1} $d v^{(1)} = W(c_{p,1}) + W(c_{p+1,0})
= d\omega + W(c_{p,1})$. Moreover, by Lemma \ref{W2} and by
$(a)_1$, we get $\|v^{(1)}\|^2 \leq cst \|c_{p,0}\|^2 \leq c_3 \|d
\omega\|^2$.

\smallskip

\textbf{Induction hypothesis}: (for $1\leq q <p+1$) there exist
positive constants $c_1$, $c_2$ and $c_3$ depending only on $n$,
$p$, $\kappa$ and $\eps$ such that the three following assertions
hold.
\begin{enumerate}
\item[$(a)_q$] There exists $c_{p+1-q,q-1} \in \Cq{q-1}{U}{\Lambda^{p+1-q}}$
such that \\ $d c_{p+1-q,q-1} = c_{p+1-(q-1),q-1}$ and
$\|c_{p+1-q,q-1}\|^2 \leq c_1 \| d \omega \|^2$.
\item[$(b)_q$] Let $c_{p+1-q,q} = (-1)^{q+1} q\cdot \delta c_{p+1-q,q-1}$.
We have $\delta c_{p+1-q,q} = 0 = dc_{p+1-q,q}$ and
$\|c_{p+1-q,q}\|^2 \leq c_2 \|d\omega\|^2$.
\item[$(c)_q$] Let $v^{(q)} = v^{(q-1)} + W(c_{p+1-q,q-1}) \in \Lambda^{p}(M)$. We
have \\ $d\omega = d v^{(q)} + (-1)^{q} W(c_{p+1-q,q})$ and
$\|v^{(q)}\|^2 \leq c_3 \|d\omega\|^2$.
\end{enumerate}

\smallskip

\textbf{Proof}: suppose the hypothesis of induction is satisfied
for some $1 \leq q \leq p$ and let us show it holds for $q+1$. By
$(b)_q$, Lemma \ref{poincare2} and Lemma \ref{lmmmu}, there exists
$c_{p-q,q} \in \Cq{q}{U}{\Lambda^{p-q}}$ and $\mu >0$ such that
$dc_{p-q,q} = c_{p+1-q,q}$ and $\|c_{p-q,q}(I)\|^2_{L^2(U_I)} \leq
\frac{2}{\mu} \|c_{p+1-q,q}(I)\|^2_{L^2(U_I)}$. Combined with
$(b)_q$ this implies that $\|c_{p-q,q}\|^2 \leq
\frac{2}{\mu}\|c_{p+1-q,q}\|^2 \leq c_1 \|d \omega \|^2$ which is
$(a)_{q+1}$. Let us consider now $$c_{p-q,q+1} = (-1)^{q} (q+1)
\delta c_{p-q,q}$$ then clearly $\delta c_{p-q,q+1} = 0$ and $d
c_{p-q,q+1} = (-1)^{q}(q+1) \delta c_{p+1-q,q} = 0$ by $(b)_q$.
Moreover, $\|c_{p-q,q+1}\|^2 \leq cst \|c_{p-q,q}\|^2 \leq c_2 \|d
\omega \|^2$ by $(a)_{q+1}$. This concludes the proof of
$(b)_{q+1}$. Finally, if $v^{(q+1)} = v^{(q)} + W(c_{p-q,q})$ we
obtain with $(c)_q$ and Lemma \ref{W1} that
\begin{eqnarray} \nonumber
    d\omega & =  & d v^{(q+1)} - d(W(c_{p-q,q})) + (-1)^{q}
    W(c_{p+1-q,q})\\
    \nonumber &=& d v^{(q+1)} - (q+1) W(\delta c_{p-q,q})
     - (-1)^q W(d c_{p-q,q}) + (-1)^{q} W(c_{p+1-q,q})\\
     &=& d v^{(q+1)} + (-1)^{q+1}W( c_{p-q,q+1}) \text{.}\nonumber
\end{eqnarray}
Finally, thanks to Lemma \ref{W2}, $(c)_q$ and $(a)_{q+1}$ we
obtain that $\|v^{(q+1)}\|^2 \leq cst (\|v^{(q)}\|^2 +
\|c_{p-q,q}\|^2) \leq c_3 \|d\omega\|^2$. This concludes the
induction.

\smallskip

\textbf{End of the induction}: (for $q=p+1$) we get $c_{0,p+1} \in
\Cq{p+1}{U}{\Lambda^{0}}$ such that $d c_{0,p+1} = 0$. This
implies in particular that $c_{0,p+1} \in i(\CqR{p+1}{\U})$.
Moreover by the proof of the de Rham Theorem seen in Section
\ref{pfdeR}, the cochain $c_{0,p+1}$ represents the same
cohomology class as $d \omega$ i.e. there exists $\gamma \in
\CqR{p}{\U}$ such that $i(\delta \gamma) = c_{0,p+1}$.

\begin{dfn} We define $\D \omega$ as the unique \Cech $p$-cochain
in $\delta^*\CqR{p+1}{\U}$ such that $i\left(\delta
\D(\omega)\right) = c_{0,p+1}$.
\end{dfn}

We prove now $(i)_{\D}$ and $(ii)_{\D}$. Firstly, by  $(b)_{p+1}$
of the induction we get that there exists a constant $c_1$
depending only on $n$, $p$, $\kappa$ and $\eps$ such that
$$\|\delta \D (\omega)\|^2 \leq cst \|c_{0,p+1}\|^2 \leq c_1 \|d
\omega\|^2$$ and this proves $(i)_{\D}$. Secondly, by $(c)_{p+1}$
we can write
\begin{eqnarray} \label{domega}
d \omega = d v^{(p+1)} + (-1)^{p+1} W(\delta \D (\omega)) = d
v^{(p+1)} + \frac{(-1)^{p+1}}{p+1}d W( \D (\omega))
\end{eqnarray}
where we used Lemma \ref{W1} and the fact that
$d(i(\D(\omega)))=0$ in the last equality. Moreover, as $\omega$
is coexact, and if $coex(\cdot)$ denotes the coexact part of a
form given by the Hodge decomposition, we deduce that
\begin{eqnarray}\nonumber
\omega =coex( v^{(p+1)}) +  \frac{(-1)^{p+1}}{p+1}coex\left( W( \D
(\omega))\right)\text{.}
\end{eqnarray}
Therefore, by Lemma \ref{W2} and using this last equality we
obtain
\begin{eqnarray}\label{domegaN}
\|\D(\omega)\| \geq cst \|W(\D(\omega))\| \geq cst (\|\omega\| -
\|v^{(p+1)}\|) .
\end{eqnarray}
Finally, by $(c)_{p+1}$ there exists $C'$ depending only on $n$,
$p$, $\kappa$ and $\eps$ such that $\|\D(\omega)\| \geq cst
(\|\omega\| - C'\|d\omega\| )$. Let then $\Lambda =
\frac{1}{4C'^2}$ so that if $\|d\omega\|^2 \leq \Lambda
\|\omega\|^2$ then $\|\D(\omega)\| \geq c_2 \|\omega\|$ which is
the requested inequality in $(ii)_{\D}$. $\square$


\subsection{From \Cech cochains to smooth forms} \label{constrS}
In this section, we are going to construct $$\S :
\delta^*\CqR{p+1}{\U} \to d^* \Lambda^{p+1}(M)$$ such that there
exist positive constants $c'_1$, $c'_2$ and $\Lambda'$ depending
only on $n$, $p$, $\kappa$ and $\eps$ such that
\begin{enumerate}
    \item[$(i)_{\S}$] $\|d \S (c) \|^2 \leq c'_1 \|\delta c\|^2$,
    for any $c \in \delta^*\CqR{p+1}{\U}$,
    \item[$(ii)_{\S}$] $\|\S c \|^2 \geq c'_2 \|c\|^2$,
    for any $c \in \delta^*\CqR{p+1}{\U}$ satisfying
    $\| \delta c \|^2 \leq \Lambda' \| c\|^2$.
\end{enumerate}
The construction of $\S$ is similar to the construction of $\D$.
The main difference is that the Whitney map is not the suitable
tool to obtain $(ii)_{\S}$. So we have to do a first induction to
construct $\S$ and a second induction (slightly different) to
prove $(ii)_{\S}$. We begin by adjusting Lemma \ref{poinc} to our
purpose.
\begin{lmm}\label{poinc2}
    Let $\U$ be a contractible cover and $\{\varphi_j\}$ a
    partition of unity subordinated to $\U$. Let $\nu >0$ such
    that $|\{j : U_j \cap U_I \neq \emptyset\}| \leq \nu$ for any
    $I \in S_{k}(\U)$ and any $k=0,\ldots ,n$. Let $c \in
    \Cq{q}{U}{\Lambda^p}$ ($q \geq 1$) such that $\delta c =0$.
    Then there exists $b \in \Cq{q-1}{U}{\Lambda^p}$ such that
    $\delta b = c$ and there exist positive constants $c_1$, $c_2$ depending only on
    $\nu$ and on a bound on $\|d\varphi_j\|_{\infty}$ such that
    \begin{enumerate}
        \item[(i)] $\|b\|^2 \leq c_1 \|c\|^2$
        \item[(ii)] $\|d b \|^2 \leq c_2 (\|c\|^2 + \|dc\|^2)$
    \end{enumerate}
\end{lmm}

\textbf{Proof}: a suitable $b$ is given by Lemma A.4.1 in
\cite{Go} and defined by $$b(I) = \sum_{j \text{ s.t. } U_j \cap
U_I \neq \emptyset} \varphi_j \cdot c(\{j\}\cup I)$$ so that $b$
verifies already $\delta b = c$. Then $(i)$ is an immediate
consequence of the definition of $b$ and $\nu$. It remains to show
$(ii)$. We have $\|db\|^2 = \sum_{I \in S_{q-1}(\U)} \|d
b(I)\|^2$. Moreover

\begin{eqnarray}
    \|d b (I)\|^2   &=& \left\|\sum_{j \text{ s.t. } U_j \cap U_I \neq \emptyset}
                    d\varphi_j\wedge c(\{j\}\cup I) + \varphi_j d c(\{j\}\cup I) \right\|^2 \nonumber \\
                    &\leq& 2 \nu \sum_{j \text{ s.t. } U_j \cap U_I \neq \emptyset}
                    \|d\varphi_j\wedge c(\{j\}\cup I)\|^2 +
                    \|\varphi_j d c(\{j\}\cup I) \|^2 \nonumber
\end{eqnarray}
and this implies the claim. $\square$

\begin{rmk}
In the sequel, we will consider a partition of unity
$\{\varphi_j\}$ subordinated to an open cover made of balls of
radius $\eps$, so that we can find a bound on
$\|d\varphi_j\|_{\infty}$ depending only on $\eps$. In particular,
this bound will be replaced by a constant depending only on
$\eps$.
\end{rmk}

\subsubsection*{Construction by induction of $ \S (\cdot )$}

Let us now proceed to the construction of $\S$ and to the proof of
$(i)_{\S}$. Let $c \in \delta^*\CqR{p+1}{\U}$. Then $\delta c$ is
an exact \Cech $(p+1)$-cochain.

\textbf{First step of induction}: define $c_{0,p+1} \in
\Cq{p+1}{U}{\Lambda^{0}}$ by $c_{0,p+1} = i(\delta c)$ i.e.
$c_{0,p+1}(I) = \delta c(I)$ for any $I \in S_{p+1}(\U)$. Clearly,
$\delta c_{0,p+1} = 0 = d c_{0,p+1}$. Then there exist positive
constants $c'_1$, $c'_2$ depending only on $n$, $p$, $\kappa$ and
$\eps$ such that
\begin{enumerate}
\item[$(a')_1$] there exists $c_{0,p} \in \Cq{p}{U}{\Lambda^{0}}$
such that $\delta c_{0,p} = c_{0,p+1}$ and $\|c_{0,p}\|^2 \leq
c'_1 \|\delta c\|^2$.
\item[$(b')_1$] Let $c_{1,p} = d c_{0,p}$. Then $\delta c_{1,p} = 0$
and $\|c_{1,p}\|^2 \leq c'_2 \|\delta c\|^2$.
\end{enumerate}
Indeed, $(a')_1$ is a direct consequence of Lemma \ref{poinc2} as
$\delta c_{0,p+1} = 0$ and of (\ref{updel}). The bound on the norm
of $dc_{0,p}$ follows also from Lemma \ref{poinc2} as $d c_{0,p+1}
= 0$. Finally, we have $\delta c_{1,p} = d \delta c_{o,p} = d
c_{0,p+1} = 0$.

\smallskip

\textbf{Induction hypothesis}: (for $1 \leq q < p+1$) there exist
positive constants $c'_1$, $c'_2$ depending only on $n$, $p$,
$\kappa$ and $\eps$ such that
\begin{enumerate}
\item[$(a')_q$] there exists $c_{q-1,p+1-q} \in \Cq{p+1-q}{U}{\Lambda^{q-1}}$
such that \\ $\delta c_{q-1,p+1-q} = c_{q-1,p+1-(q-1)}$ and
$\|c_{q-1,p+1-q}\|^2 \leq c'_1 \|\delta c\|^2$.
\item[$(b')_q$] Let $c_{q,p+1-q} = d c_{q-1,p+1-q}$. Then $\delta c_{q,p+1-q} = 0$
and $\|c_{q,p+1-q}\|^2 \leq c'_2 \|\delta c\|^2$.
\end{enumerate}

\smallskip

\textbf{Proof}: suppose the hypothesis of induction is verified
for some $1\leq q  \leq p$ and let us show it holds for $q+1$. By
$(b')_q$ and by Lemma \ref{poinc2} there exists $c_{q,p-q} \in
\Cq{p-q}{U}{\Lambda^{q}}$ such that $\delta c_{q,p-q} =
c_{q,p+1-q}$ and $\|c_{q,p-q}\|^2 \leq cst \|c_{q,p+1-q}\|^2$.
Combined with $(b')_q$, this implies $(a)_{q+1}$. Moreover, let us
consider $c_{q+1,p-q}=d c_{q,p-q}$. Then, by definition of
$c_{q,p+1-q}$ we have $\delta c_{q+1,p-q} = d \delta c_{q,p-q} = d
c_{q, p+1-q} = 0$. Finally, by Lemma \ref{poinc2}, we have
$\|c_{q+1,p-q}\|^2 \leq cst (\|c_{q,p+1-q}\|^2 + \|d
c_{q,p+1-q}\|^2)$. As we have $d c_{q,p+1-q} = 0$ and by $(b')_q$,
we get $\|c_{q+1,p-q}\|^2 \leq c'_2 \|\delta c\|^2$. This
concludes the induction.

\smallskip

\textbf{End of the induction}: (for $q = p+1$) we obtain
$c_{p+1,0} \in \Cq{0}{U}{\Lambda^{p+1}}$ such that $\delta
c_{p+1,0} = 0$. This implies that $c_{p+1,0}$ is the restriction
of a well-defined $(p+1)$-form and by the de Rham Theorem as
$\delta c$ is exact, the $0$-cochain $c_{p+1,0}$ is exact and is
the restriction of an exact $(p+1)$-form.

\begin{dfn}
    Let $\S (c)\in d^* \Lambda^{p+1}(M)$ be the unique coexact $p$-form such that
    $r(d \S(c)) = c_{p+1,0}$.
\end{dfn}

An immediate consequence of the induction is $(i)_{\S}$. Indeed,
from $(b')_{p+1}$ and Remark \ref{norm} follows that there exists
a positive constant $c'_1$ depending only on $n$, $p$, $\kappa$
and $\eps$ such that $\|d \S(c)\|^2 \leq c'_1 \|\delta c\|^2$.

Let us now proceed to a second induction in order to prove
$(ii)_{\S}$. The goal is to construct $b \in \CqR{p}{\U}$ such
that $\delta b = \pm \delta c$ and $\|b\| \leq cst (\|\S(c)\| +
\|\delta c\| )$ where $cst$ is a positive constant depending only
on $n$, $p$, $\kappa$ and $\eps$. These are in fact the
corresponding equations for (\ref{domega}) and (\ref{domegaN}) in
the discretizing part. In the induction, we will use the $c_{r,s}$
appearing in the construction of $\S$.

\textbf{First step of induction}: define $b_{p,0} = r(\S(c)) -
c_{p,0} \in \Cq{0}{U}{\Lambda^{p}}$. We have $d b_{p,0} =
c_{p+1,0} - d c_{p,0} = 0$. Then there exist positive constants
$c''_1$, $c''_2$ depending only on $n$, $p$, $\kappa$ and $\eps$
such that
\begin{enumerate}
\item[$(a'')_1$] there exists $b_{p-1,0} \in
\Cq{0}{U}{\Lambda^{p-1}}$ such that $d b_{p-1,0} = b_{p,0}$ and
$\|b_{p-1,0}\|^2 \leq c''_1 (\|\S(c)\| + \|\delta c\|)$.
\item[$(b'')_1$] Let $b_{p-1,1} = \delta b_{p-1,0} + c_{p-1,1}$.
Then we have $d b_{p-1,1} = 0$ and $\|b_{p-1,1}\| \leq
c''_2(\|\S(c)\| + \|\delta c\| )$.
\end{enumerate}
Indeed, as $p \geq 1$ and $ d b_{p,0} =0$, by Lemma
\ref{poincare2} there exists $b_{p-1,0} \in
\Cq{0}{U}{\Lambda^{p-1}}$ such that $d b_{p-1,0} = b_{p,0}$ and
$\|b_{p-1,0}\| \leq cst \|b_{p,0}\|$. By definition of $b_{p,0}$
and by $(a')_{p+1}$ of the previous induction we obtain then
$(a'')_1$. Let us consider now $b_{p-1,1} = \delta b_{p-1,0} +
c_{p-1,1}$. Then we have $d b_{p-1,1} = \delta b_{p,0} + c_{p,1} =
- \delta c_{p,0} + c_{p,1}  = 0$. Finally, by construction and by
(\ref{updel}) $\|b_{p-1,1}\| \leq cst (\|b_{p-1,0}\| +
\|c_{p-1,1}\|)$. This last inequality combined with $(a'')_1$ and
$(a')_p$ leads to $(b'')_1$.

\smallskip

\textbf{Induction hypothesis}: (for $1 \leq q < p-1$) there exist
positive constants $c''_1$, $c''_2$ depending only on $n$, $p$,
$\kappa$ and $\eps$ such that
\begin{enumerate}
\item[$(a'')_q$] there exists $b_{p-q,q-1} \in
\Cq{q-1}{U}{\Lambda^{p-q}}$ such that $d b_{p-q,q-1} =
b_{p-(q-1),q-1}$ and $\|b_{p-q,q-1}\|^2 \leq c''_1 (\|\S(c)\| +
\|\delta c\|)$.
\item[$(b'')_q$] Let $b_{p-q,q} = \delta b_{p-q,q-1} + (-1)^{q+1} c_{p-q,q}$.
Then we have $d b_{p-q,q} = 0$ and $\|b_{p-q,q}\| \leq
c''_2(\|\S(c)\| + \|\delta c\| )$.
\end{enumerate}

\smallskip

\textbf{Proof}: suppose the induction hypothesis holds for some $
1 \leq q \leq p-1$ and let us show it holds for $q+1$. By
$(b'')_q$ and Lemma \ref{poincare2} there exists $ b_{p-(q+1),q}
\in \Cq{q}{U}{\Lambda^{p-(q+1)}}$ such that $d b_{p-(q+1),q} =
b_{p-q,q}$ and $\|b_{p-(q+1),q}\|^2 \leq cst \|b_{p-q,q}\|^2$ and
it suffices to use $(b'')_q$ to obtain $(a'')_{q+1}$. Then
consider $b_{p-(q+1),q+1} = \delta b_{p-(q+1),q} +
(-1)^qc_{p-(q+1),q+1}$. We have
\begin{eqnarray}
d b_{p-(q+1),q+1}   &=& \delta  b_{p-q,q} +
                        (-1)^q c_{p-q,q+1} \nonumber \\
                    &=& \delta (\delta b_{p-q,q-1}+ (-1)^{q+1}
                    c_{p-q,q}) + (-1)^q \delta c_{p-q,q} \nonumber\\
                    &=& 0 . \nonumber
\end{eqnarray}
Finally, by construction of $b_{p-(q+1),q+1}$ we have
$$\|b_{p-(q+1),q+1}\| \leq cst (\|b_{p-(q+1),q}\| +
\|c_{p-(q+1),q+1}\|)$$ and with $(a'')_{q+1}$ and $(a')_{p-q}$ we
obtain $(b'')_{q+1}$. This ends the induction.

\smallskip

\textbf{End of the induction}: (for $q=p$) we obtain $b_{0,p} \in
\Cq{p}{U}{\Lambda^0}$ such that $d b_{0,p} = 0$ i.e. $b_{0,p} \in
\CqR{p}{\U}$ and $\delta b_{0,p} = (-1)^{p+1} \delta c_{0,p} =
(-1)^{p+1} c_{0,p+1} = (-1)^{p+1} \delta c$. Hence, $b_{0,p}$ and
$c$ have same coexact part and as $c$ is already coexact we obtain
by $(b'')_{p}$, $\|c\| \leq \|b_{0,p}\| \leq cst (\|\S(c)\| +
\|\delta c\|)$. In particular,
\begin{eqnarray} \nonumber
\|\S(c)\| \geq \frac{1}{cst}\|c\| - \|\delta c\|
\end{eqnarray}
then let $\Lambda' = \frac{1}{4cst^2}$ so that if $\|\delta c\|^2
\leq \Lambda \|c\|^2$ then $\|\S(c)\| \geq c'_2 \|c\|$. This ends
the proof of $(ii)_{\S}$. $\square$


\subsection{Upper bounds on the spectra} \label{secup}

\begin{lmm} \label{upper1} Let $(M^n,g)$ be a compact connected Riemannian manifold
and let $\U$ be a finite contractible open cover of $M$ such that
there exists $\nu >0$ such that $| \{j : U_j \cap U_{I} \neq
\emptyset \}| \leq \nu$ for any $I \in S_{q}(\U)$ and any $ q \geq
0$. Then there exists a positive constant $c$ depending only on
$\nu$ and $p$ such that $\lambda_{k,q}(\U) \leq c $ for any $k =
1, \ldots, |S_q(\U)| - \check{b}_q(\U)$.
\end{lmm}
\textbf{Proof}: it suffices to show the result for the spectrum of
$\delta^* \delta$ on $\delta^* \CqR{p+1}{\U}$. We are going to
show that there exists a positive constant depending only on $\nu$
and $p$ such that for any $b \in  \CqR{p}{\U}$
\begin{eqnarray} \label{updel}
\| \delta b\|^2 \leq cst \|b\|^2
\end{eqnarray}
and then the variational characterization of the spectrum of
$\delta^* \delta$ will imply the claim. Recall that $\delta b(I) =
\sum_{i\in I} \epsilon(i,I\setminus i) b(I\setminus i)$ where
$\epsilon(i,I\setminus i)$ denotes the signature of the
permutation ordering $\{i\} \cup (I\setminus i)$ to obtain $I$ and
$I \in S_{p+1}(\U)$. Hence
\begin{eqnarray}\nonumber
|\delta b(I)|^2 \leq (p+2) \sum_{i\in I} |b(I\setminus i)|^2 .
\end{eqnarray}
This implies that
\begin{multline}\nonumber
    \|\delta b\|^2  = \sum_{I \in S_{p+1}(\U)}|\delta b(I)|^2
                    \leq (p+2) \sum_{I \in S_{p+1}(\U)}
                    \sum_{i\in I} |b(I\setminus i)|^2 \\
                    \leq (p+2) \nu \sum_{J \in S_{p}(\U)}|b(J)|^2
                    = (p+2) \nu \|b\|^2
\end{multline}
which is the claim. $\square$

\begin{lmm}\label{upper2} Let $(M,g) \in \K$ and $X$ an $\eps$-discretization
with $0 < \eps \leq r_0$. Let $1 \leq p \leq n-1$. Then there
exists a positive constant $c'$ depending only on $n$, $p$,
$\kappa$ and $\eps$ such that $\lambda_{k,p}(M) \leq c'$ for any
$k \leq |S_{p}(\U_X)|-\check{b}_p(\U_X)$.
\end{lmm}
\textbf{Proof}: it suffices to show the result for $k =
|S_{p}(\U_X)|-\check{b}_p(\U_X)$. By a theorem  of Abresch (see
\cite{CFG}, Theorem 1.12) there exists a Riemannian metric
$\tilde{g}$ on $M$ such that
\begin{enumerate}
\item[(a)] $e^{-\frac{1}{4}} g \leq \tilde{g} \leq e^{\frac{1}{4}} g$
\item[(b)] $|\nabla^g - \nabla^{\tilde{g}}| \leq \frac{1}{4}$
\item[(c)] $|K_{\tilde{g}}| \leq \tilde{\kappa}(n,\kappa) \text{
and }|\nabla^{\tilde{g}} R_{\tilde{g}}| \leq K(n,\kappa)$
\end{enumerate}
where $\tilde{\kappa}$ and $K$ depend only on $n$ and $\kappa$. By
Proposition 3.3. of \cite{Do}, there exist a positive constant $c$
depending only on $e^{\frac{1}{4}}$ such that $$\lambda_{k,p}(M,g)
\leq c  \lambda_{k,p}(M,\tilde{g}). $$ Therefore it suffices to
show the claim for $(M,\tilde{g})$. By Remark \ref{cardX} and by
construction of $\tilde{g}$, there exists a positive constant $d$
depending only on $n$, $p$, $\kappa$, $\eps$ such that
$|S_{p}(\U_X)| \leq d Vol(M,\tilde{g})$. Moreover, there exist
$\alpha > 0$ depending only on $p$, $n$, $\kappa$ and $\eps$ such
that if $Y$ is an $\alpha$-discretization of $(M,\tilde{g})$ then
$|Y| \geq |S_{p}(\U_X)|$ and $\check{b}_p(\U_Y) =
\check{b}_p(\U_X)$. Consider then the disjoint balls (for
$\tilde{g}$) centered at $y \in Y$ of radius $\frac{\alpha}{2}$.
From Proposition 2.3. of \cite{Do}, on any of these balls there
exists a $p$-form $\omega_y$ which is zero on the boundary of the
ball, so that we can extend $\omega_y$ by zero to obtain a
$p$-form on $M$ also denoted $\omega_y$ such that
\begin{eqnarray}\label{eqnmu}
\frac{\|d \omega_y\|^2_{\tilde{g}} + \|d^*_{\tilde{g}} \omega_y
\|^2_{\tilde{g}}}{\|\omega_y\|^2_{\tilde{g}}} \leq
\mu(n,p,\kappa,\eps) \end{eqnarray} where $\mu(n,p,\kappa,\eps)$
is a positive constant depending only on $n$, $p$, $\kappa$ and
$\eps$. Moreover, we can choose $\omega_y$ such that $\|\omega_y\|
= 1$.

\smallskip

Let then $V$ the vector subspace of $p$-forms spanned by
$\{\omega_y : y \in Y\}$. By construction, $\omega_y$ is
orthogonal to $\omega_x$ if $x \neq y$. In particular, $V$ is of
dimension $|Y|$. Therefore, by the variational characterization of
the spectrum, we obtain
\begin{eqnarray} \label{eqnl}
\lambda_{|Y| - \check{b}_p(\U_Y) ,p}(M,\tilde{g}) \leq
\max\left\{\frac{\|d \omega\|^2_{\tilde{g}} + \|d^*_{\tilde{g}}
\omega \|^2_{\tilde{g}}}{\|\omega\|^2_{\tilde{g}}} : \omega \in V
\setminus \{0\} \right\} \text{.}
\end{eqnarray}
Furthermore, if $\omega = \sum_{y\in Y} a_y \omega_y$, then as the
balls centered on $Y$ of radius $\frac{\alpha}{2}$ are disjoint
$\|\omega\|_{\tilde{g}}^2 \geq \sum_{y\in Y} a_y^2$ and combined
with (\ref{eqnmu}) this implies that
\begin{eqnarray}\label{star}
\|d\omega\|^2_{\tilde{g}} \leq \sum_{y\in Y} a_y^2 \|d
\omega_y\|^2_{\tilde{g}} \leq \mu \|\omega\|^2_{\tilde{g}}
\end{eqnarray}
and
\begin{eqnarray}\label{starstar2}
\|d^*_{\tilde{g}} \omega\|^2_{\tilde{g}} \leq \sum_{y\in Y} a_y^2
\|d^*_{\tilde{g}} \omega_y\|^2_{\tilde{g}} \leq \mu
\|\omega\|^2_{\tilde{g}}\text{.}
\end{eqnarray}

It suffices then to introduce (\ref{star}) and (\ref{starstar2})
in (\ref{eqnl}) to obtain that $$\lambda_{|Y| -
\check{b}_p(\U_Y),p}(M,\tilde{g}) \leq 2 \mu$$ and in particular
that $\lambda_{k,p}(M,g) \leq 2 c \mu$, for $k \leq
|S_{p}(\U_X)|-\check{b}_p(\U_X)$. $\square$


\subsection{Proof of the main result}

We prove now Theorem \ref{thmmain}. We will only proceed to the
proof of the inequality $\lambda_{k,p}(M) \leq c_2
\lambda_{k,p}(X)$ as the other inequality can be proved in the
same way using the corresponding results. Recall it suffices to
prove the result for $d^*d$ on coexact forms and for $\delta^*
\delta$ on coexact \Cech cochains. We proceed in two steps. Let
$\Lambda'$ given by $(ii)_{\S}$.

\smallskip

\textbf{First step}: assume $ \lambda^{\delta^*}_{k,p}(X) \geq
\Lambda'$. Then, $\lambda^{d^*}_{k,p}(M) \leq
\Lambda'^{-1}\lambda^{\delta^*}_{k,p}(X)\lambda^{d^*}_{k,p}(M)$
and by Lemma \ref{upper2} we obtain  $\lambda^{d^*}_{k,p}(M) \leq
\Lambda'^{-1} c' \lambda^{\delta^*}_{k,p}(X)$ which is the claim.

\smallskip

\textbf{Second step}: assume now $\lambda^{\delta^*}_{k,p}(X) \leq
\Lambda'$. Let us consider $c_1,\ldots,c_k \in \delta^*
\Cq{p+1}{\U}$ the \Cech
$\lambda^{\delta^*}_{1,p}(X),\ldots,\lambda^{\delta^*}_{k,p}(X)$-eigencochains
such that $(c_i,c_j)=\delta_{ij}$. Denote by $V^{k}$ the
$k$-dimensional vector subspace of $\delta^* \CqR{p+1}{\U}$ they
span. By the variational characterization of the spectrum we have
\begin{eqnarray}\nonumber
    \lambda^{\delta^*}_{k,p}(X) =
    \max\left\{\frac{\|\delta c\|^2}{\|c\|^2} : c \in
    V^k\setminus\{0\}\right\} .
\end{eqnarray}
Let us consider now $\S V^k$ the vector subspace of
$d^*\Lambda^{p+1}(M)$ spanned by $\{\S(c_1),\ldots,\S(c_k)\}$.
Then if $\S(c) \in \S V^k$, $\S(c) = \sum_{i=1}^{k}a_i \S(c_i)$
with $c = \sum_{i=1}^{k}a_i c_i \in V^k$. So that we have
$\|\delta c \|^2 \leq \lambda_{k,p}^{\delta^*}(X) \|c\|^2 \leq
\Lambda' \|c\|^2$. Therefore, by $(ii)_{\S}$ we obtain
\begin{eqnarray}\label{eqn}
\|\S(c)\|^2 \geq c'_2 \|c\|^2
\end{eqnarray}
and this says in particular that $\S V^k$ is of dimension $k$.
Using the variational characterization of $\lambda_{k,p}^{d^*}(M)$
we get
\begin{multline}\nonumber
    \lambda_{k,p}^{d^*}(M)  \leq
    \max \left\{ \frac{\|d \omega\|^2}{\|\omega\|^2} : \omega \in \S V^k\setminus\{0\} \right\}
     \\ =  \max \left\{ \frac{\|d \S(c)\|^2}{\|\S(c)\|^2} : c \in V^k\setminus\{0\}
    \right\} \text{.}
\end{multline}
Finally, (\ref{eqn}) and $(i)_{\S}$ imply that $\frac{\|d
\S(c)\|^2}{\|\S(c)\|^2} \leq \frac{c'_1}{c'_2} \frac{\| \delta
c\|^2}{\|c\|^2}$ so that we obtain
\begin{eqnarray}
    \lambda_{k,p}^{d^*}(M)  \leq
    \frac{c'_1}{c'_2} \max \left\{ \frac{\|\delta c \|^2}{\|c\|^2} : c \in V^k\setminus\{0\}
    \right\} = \frac{c'_1}{c'_2} \lambda^{\delta^*}_{k,p}(X)
\end{eqnarray}
which concludes the proof. $\square$


%% file: DiscFormsApplication.tex
\section{Applications} \label{sec4}

In this section, we develop several consequences of Theorem
\ref{thmmain} or of the methods used to prove Theorem
\ref{thmmain}.

\subsection{A lower bound for the spectrum of the Laplacian on differential
forms} \label{sec41}

The goal of this section is to prove the following theorem.

\begin{thm}\label{lambda1}
    Let $(M,g) \in \K$. Let $1 \leq p \leq n-1$. Then there exists
    a positive constant $c(n,p,\kappa,r_0)$ depending only on $n$, $p$, $\kappa$ and
    $r_0$ such that $$\lambda_{1,p}(M) \geq \frac{c(n,p,\kappa,r_0)}{Vol(M)
    e^{Vol(M)}}$$
    where $Vol(M)$ denotes the volume of $(M,g)$.
\end{thm}

By Theorem \ref{thmmain}, it suffices to choose a suitable
discretization $X$ of $M$ and prove then a similar result for
$\lambda_{1,p}(X)$. To that aim we need the following lemma.

\begin{lmm} \label{lmmT}
    Let $A : \R^m \to \R^n$ be a linear operator with matrix
    coefficients (in the canonical bases) in $\{-1,0,1\}$. Suppose
    there exists an integer $k$ such that any column and any row
    has at most $k$ non-zero coefficients. Then, there exists $B :
    \R^n \to \R^m$ such that $ABAv = Av$ for any $v \in \R^m$ and
    $$\|Bu\|^2 \leq n k^{2n} \|u\|^2$$ for any $u \in \R^n$.
\end{lmm}

\begin{rmk} \label{faute} In \cite{Tr}, the author proves a similar result (see Lemma
A.5 in \cite{Tr}) but with a better constant for the matrix norm
of $B$. He asserts that $\|B u\|^2 \leq  c(k) m \|u\|^2$. The
following example shows that the constant in Trèves' result is not
suitable. Consider the matrix $A$ with $m$ columns and $m-1$ rows
given by
\begin{equation*}
    A =
    \begin{pmatrix}
    1 & -1 & 0  & \ldots & 0 \\
    0 & 1 & -1  & \ldots & 0 \\
    \vdots & \ddots & \ddots & \ddots &  \vdots \\
    0 & \ldots & 0 & 1 & -1  \\
    \end{pmatrix}
\end{equation*}
and consider $v = \sum_{i=1}^{m} i e_i$ in $\R^m$. Then $Av =
-\sum_{i=1}^{m-1}e_i$ in $\R^{m-1}$. So that $\|Av\|^2 = m-1$. An
easy calculation shows that if we choose the $m-1$ first columns
of $A$ to span $Im(A)$ then $BAv = \sum_{i=1}^{m-1} -(m-i)e_i $ in
$\R^m$. Hence $\|BAv\|^2 = \frac{(m-1)m(2m-1)}{6} =
\frac{m(2m-1)}{6}\|Av\|^2$ which contradicts Lemma A.5 in
\cite{Tr} (here $k=2$). The assertion $A.44$ in \cite{Tr} is
wrong. It is not clear to us how we can correct this mistake. We
think that we should replace $k^{2n}$ by $n^{l}$ for a suitable
$l$ in Lemma \ref{lmmT} but we cannot prove it yet.
\end{rmk}

\textbf{Proof of Lemma \ref{lmmT}}: let $r$ be the dimension of
$Im(A)$. Without lost of generality we can suppose that the $r$
first columns $\{a_1, \ldots a_r \}$ of $A$ span $Im(A)$. Then
define $B$ as follows. On the orthogonal complement of $Im(A)$ let
$B=0$. Moreover, if $u = Av$ then write $u$ in the basis
$\{a_1,\ldots,a_r\}$ of $Im(A)$, $u=\sum_{i=1}^r u_i a_i$ and
define $Bu = \sum_{i=1}^r u_i e_i$ where $\{e_i\}$ denotes the
canonical basis of $\R^m$. An immediate consequence of the
definition of $B$ is that $ABAv = Av$. Moreover, $\|Bu\|^2 =
\sum_{i=1}^r u_i^2$. Let us show now that
\begin{eqnarray} \label{uicarre}
    u_i^2 \leq k^{2n} \|u\|^2 .
\end{eqnarray}
This will imply $\|Bu\|^2 \leq r k^{2n} \|u\|^2 \leq n k^{2n}
\|u\|^2$ which is the claim.

We prove (\ref{uicarre}) for $i=1$. Let $V_1$ the vector space
spanned by $\{a_2,\ldots,a_r\}$ and let $V_1^{\perp}$ its
orthogonal complement in $Im(A)$. Consider $P_1 : Im(A) \to
V_1^{\perp}$ the orthogonal projection onto $V_1^{\perp}$. We have
$P_1(u) = u_1 P_1(a_1)$ so that
\begin{eqnarray} \label{u1carre}
    u_1^2 = \frac{\|P_1(u)\|^2}{\|P_1(a_1)\|^2} \leq
    \frac{\|u\|^2}{\|P_1(a_1)\|^2} .
\end{eqnarray}
We can write $P_1(a_1) = a_1 + \alpha_2 a_2 + \ldots + \alpha_r
a_r$ with $\left(P_1(a_1)|a_j\right) = 0$ for $j=2,\ldots,r$ and
$\left(P_1(a_1)|a_1\right) = \|P_1(a_1)\|^2$. In matrix form we
obtain
\begin{equation*}
    \begin{pmatrix}
    \|a_1\|^2 & (a_1|a_2) & \ldots & (a_1|a_r)\\
    (a_1|a_2) & \|a_2\|^2 & \ldots & (a_2|a_r)\\
    \vdots    & \vdots    & \ddots & \vdots   \\
    (a_1|a_r) & (a_2|a_r) & \ldots & \|a_r\|^2\\
    \end{pmatrix}
    \begin{pmatrix}
    1 \\
    \alpha_2 \\
    \vdots \\
    \alpha_r\\
    \end{pmatrix}
    =
    \begin{pmatrix}
    \|P_1(a_1)\|^2 \\
    0 \\
    \vdots \\
    0 \\
    \end{pmatrix}
\end{equation*}
and if we call $P$ the matrix $r\times r$ above and $Q$ the
submatrix of $P$ obtained by removing the first row and the first
column of $P$ we get that
\begin{eqnarray}\nonumber
    \|P_1(a_1)\|^2 = \frac{|\det(P)|}{|\det(Q)|} .
\end{eqnarray}
As $\{a_1,\ldots a_r\}$ are linearly independent, $\det(P) \neq
0$. Moreover, $P$ is a matrix with integer coefficients so that
$|\det(P)| \geq 1$. It remains to find an upper bound for
$|\det(Q)|$. So, we are going to prove by induction that the
minors of $P$ of size $l\times l$ are bounded above by $k^{2l-1}$.

\smallskip

The first step of induction asserts that the minors of $P$ of size
$1 \times 1$ are bounded above by $k$. This is a direct
consequence of the assumption that each column of $A$ has at most
$k$ non-zero coefficients. Suppose then that the minors of $P$ of
size $l\times l$ are bounded above by $k^{2l-1}$. Consider then
$D$ a minor of $P$ of size $(l+1) \times (l+1)$. Then $D$ can be
written as $$D = \sum_{j=1}^{l+1} c_j D_j$$ where
$(c_1,\ldots,c_{l+1})$ is a part of a line of $P$ and $D_j$ is a
minor of $P$ of size $l\times l$. By construction of $P$, the
coefficients $c_j$ can be written as follows. There exists $1 \leq
J \leq r$ such that $$c_j = (a_J|a_{i_j}) \text{ for a suitable
}i_j$$ so that
\begin{eqnarray*}
    |D| &=& \left|\sum_{j=1}^{l+1}(a_J|a_{i_j})D_j \right|
        = \left| \sum_{i=1}^{n}(a_J|e_i)\sum_{j=1}^{l+1}(e_i|a_{i_j})D_j \right| .\\
\end{eqnarray*}
But by assumption, the $i^{\text{th}}$ row of $A$ has at most $k$
coefficients of absolute value $1$ and by induction hypothesis we
get $|\sum_{j=1}^{l+1}(e_i|a_{i_j})D_j | \leq k\cdot k^{2l-1}$.
Moreover, by assumption the $J^{\text{th}}$ column of $A$ has at
most $k$ coefficients of absolute value $1$ and with the previous
remark this implies
\begin{eqnarray*}
    |D| \leq k \cdot k \cdot k^{2l-1}
\end{eqnarray*}
and this ends the induction. We apply then the result to
$|\det(Q)|$ and we obtain $|\det(Q)| \leq k^{2r-3} \leq k^{2n}$.
Finally, we deduce that
\begin{eqnarray}\nonumber
\|P_1(a_1)\|^2 \geq \frac{1}{k^{2n}}
\end{eqnarray}
and combined with (\ref{u1carre}) this implies (\ref{uicarre}).
$\square$

\begin{thm}\label{lambda2}
Let $\U$ be a finite open cover of $M$ compact. Let $p \geq 0$.
Assume there exists $\nu$ such that $|\{j : U_j \cap U_I \neq
\emptyset\}| \leq \nu$ for any $I \in S_q(\U)$ and $q \geq 0$.
Then there exists a positive constant $c(\nu,p)$ depending only on
$\nu$ and $p$ such that $$\lambda_{1,p}(\U) \geq
\frac{c(\nu,p)}{|\U| \cdot  e^{|\U|}} \ .$$
\end{thm}

\textbf{Proof}: it suffices to prove the result for
$\lambda_{1,p}^{\delta^*}(\U)$. By the variational
characterization of the spectrum, we have
\begin{eqnarray}\nonumber
    \lambda_{1,p}^{\delta^*}(\U) =
    \min_V \max\left\{\frac{\|\delta c\|^2}{\| c\|^2} : c \in V\setminus \{0\} \right\}
\end{eqnarray}
where $V$ ranges over all $1$-dimensional vector subspaces of
$\delta^* \CqR{p+1}{\U}$. As in Proposition 3.1 of \cite{Do}, we
can get from the above characterization the following description
\begin{eqnarray}\nonumber
    \lambda_{1,p}^{\delta^*}(\U) =
    \min_V \max\left\{\frac{\|\delta c\|^2}{\| b\|^2} : \delta b = \delta c \ , \text{ and }\delta c \in V \right\}
\end{eqnarray}
where $V$ ranges over all $1$-dimensional vector subspaces of
$\delta \CqR{p}{\U}$. In particular, if we consider $V$ that
realizes the minimum, then
\begin{eqnarray}\label{charL}
    \lambda_{1,p}^{\delta^*}(\U) =
    \max\left\{\frac{\|\delta c\|^2}{\| b\|^2} : \delta b = \delta c \
    , \text{ and }\delta c \in V \right\} .
\end{eqnarray}

Consider then the canonical basis of $\CqR{q}{\U}$ given by
\begin{eqnarray}\nonumber
\{ e_I : S_q(\U) \to \R , I \in S_q(\U) \text{ such that } e_I(J)
= \delta_{IJ} \} .
\end{eqnarray}
In this bases, the matrix of $\delta : \CqR{p}{\U} \to
\CqR{p+1}{\U}$ has coefficients in $\{-1,0,1\}$ and has at most
$K(\nu,p) = \max\{\nu,p+2\}$ non-zero coefficients by row and by
column. Hence we can apply Lemma \ref{lmmT} to $\delta$ to obtain
that for any $c \in \CqR{p}{\U}$, there exists $b \in \CqR{p}{\U}$
such that $\delta b = \delta c$ and
\begin{eqnarray} \label{bdc}
\|b\|^2 \leq |S_{p+1}(\U)|K(\nu,p)^{|S_{p+1}(\U)|} \|\delta c\|^2
.
\end{eqnarray}
Finally, if we introduce (\ref{bdc}) in (\ref{charL}) and by
Remark \ref{cardX}, we obtain
\begin{eqnarray}\nonumber
    \lambda_{1,p}^{\delta^*}(\U) \geq
    \frac{1}{|S_{p+1}(\U)|K(\nu,p)^{|S_{p+1}(\U)|}}
    \geq \frac{c(\nu,p)}{|\U| \cdot  e^{|\U|}} \ .\ \  \square
\end{eqnarray}

\textbf{Proof of Theorem \ref{lambda1}}: let $(M,g) \in \K$ and
$X$ a $\frac{\rho_0}{4}$-discretiza\-tion of $M$ (where $\rho_0$
is given by Lemma \ref{lmmrho}). By Theorem \ref{thmmain}, there
exists $c_1(n,p,\kappa,r_0) > 0$ such that
\begin{eqnarray} \label{eqn1}
    \lambda_{1,p}(M,g) \geq c_1 \lambda_{1,p}(X) .
\end{eqnarray}
Moreover, by Theorem \ref{lambda2} there exists
$c_2(n,p,\kappa,r_0) > 0$ such that
\begin{eqnarray} \label{eqn256}
    \lambda_{1,p}(X) \geq
    \frac{c_2}{|\U| \cdot e^{|\U|}} .
\end{eqnarray}
Finally, by Remark \ref{cardX} there exists $c_3(n,p,\kappa,r_0) >
0$ such that
\begin{eqnarray} \label{eqn3}
    |\U| \leq c_3 Vol(M) .
\end{eqnarray}
To conclude, put (\ref{eqn1}), (\ref{eqn256}) and (\ref{eqn3})
together to obtain that there exists $c(n,p,\kappa,r_0) > 0$ such
that
\begin{eqnarray} \nonumber
      \lambda_{1,p}(M,g) \geq \frac{c}{Vol(M) e^{Vol(M)}}
\end{eqnarray}
and this ends the proof. $\square$

\subsection{Whitney forms: a natural way of smoothing}\label{secW}

As suggested in \cite{DP}, a candidate for the smoothing operator
should be given by Whitney forms in the following way. Let
\begin{eqnarray} \nonumber
    \tilde{\S} : \delta^*\CqR{p+1}{\U} \to d^*\Lambda^{p+1}(M)\ ,
    \ c \mapsto \tilde{\S}(c) = coex(W(c))
\end{eqnarray}
where $W$ is the Whitney map (see Appendix \ref{whitney}). The
results of Dodziuk and Patodi in \cite{DP} concerning Whitney
forms can not be used in our context as their approximations
(obtained thanks to the heat kernel) involve the manifold itself.
More precisely, the constants there depend on the volume of the
manifold.

\smallskip

Here, we show that there exist positive constants $\tilde{c}_1$,
$\tilde{c}_2$ and $\tilde{\Lambda}$ depending only on $n$, $p$,
$\kappa$ and $\eps$ such that
\begin{enumerate}
    \item[$(i)_{\tilde{\S}}$] $\|d \tilde{\S} (c) \|^2 \leq \tilde{c}_1 \|\delta c\|^2$,
    for any $c \in \delta^*\CqR{p+1}{\U}$,
    \item[$(ii)_{\tilde{\S}}$] $\|\tilde{\S} c \|^2 \geq \tilde{c}_2 \|c\|^2$,
    for any $c \in \delta^*\CqR{p+1}{\U}$ satisfying
    $\| \delta c \|^2 \leq \tilde{\Lambda} \| c\|^2$.
\end{enumerate}

The inequality $(i)_{\tilde{\S}}$ is a direct consequence of Lemma
\ref{W1} and Lemma \ref{W2}. Indeed, as $dc = 0$ we have
\begin{eqnarray}\nonumber
    d \tilde{\S}(c) = d W(c) = (p+1) W(\delta c)
\end{eqnarray}
and Lemma \ref{W2} leads to $(i)_{\tilde{\S}}$.

\smallskip

The second inequality is less obvious and it can be shown adding a
point to the first induction in the construction of $\S$ in
Section \ref{constrS}. The idea is to construct a $p$-form
$u^{(0)}$ linking $\S(c)$ and $\tilde{\S}(c)$ playing the same
role as $v^{(p)}$ in the construction of $\D$ (see Section
\ref{constrD}). Then the control on the norm of $\S(c)$ (see
$(ii)_{\S}$) and a control on the norm of $u^{(0)}$ will imply the
desired inequality.

\smallskip

\textbf{Proof of} {\boldmath$(ii)_{\tilde{\S}}$}: in the "first
step of induction" (of Section \ref{constrS}), add
\begin{enumerate}
\item[$(c')_1$] there exists a positive constant $c'_3$ depending
only on $n$, $p$, $\kappa$ and $\eps$ such that if $$u^{(p)} =
(-1)^{p+2}\frac{1}{p+1}W(c_{0,p})$$ then $\|u^{(p)}\|^2 \leq c'_3
\|\delta c\|^2$ and $$d W(c) = (-1)^{p+2}(p+1) \left(d u^{(p)} +
\frac{(-1)^{(p+2)(p+1)}}{p+1} W(c_{1,p}) \right) .$$
\end{enumerate}
Indeed, by Lemma \ref{W2} and $(a')_1$, $\|u^{(p)}\|^2 \leq
cst\|c_{0,p}\|^2 \leq c'_3 \|\delta c\|^2$. Moreover, by Lemma
\ref{W1} and $(b')_1$
\begin{eqnarray} \nonumber
    d u^{(p)} & = &
    \frac{(-1)^{p+2}}{p+1} \left((p+1)W(c_{0,p+1}) + (-1)^p W(c_{1,p})
    \right)\\
    &=&
    \frac{(-1)^{p+2}}{p+1} \left(d W(c) - (-1)^{p+1} W(c_{1,p})
    \right)  .\nonumber
\end{eqnarray}
The induction hypothesis gets
\begin{enumerate}
    \item[$(c')_q$] there exists a positive constant $c'_3$ depending
    only on $n$, $p$, $\kappa$ and $\eps$ such that if
    $$u^{(p+1-q)} = u^{(p+1-(q-1))} + \frac{(-1)^{p+2}(-1)^{p+1}\ldots(-1)^{p+2-(q-1)}}{(p+1)p(p-1)\ldots(p+2-q)}W(c_{q-1,p+1-q})$$ then
    $\|u^{(p+1-q)}\|^2 \leq c'_3 \|\delta c\|^2$ and $$\frac{(-1)^{p+2}}{p+1}d W(c) =
     d u^{(p+1-q)} + \frac{(-1)^{p+2}\ldots(-1)^{p+2-q}}{(p+1)\ldots(p+2-q)}
    W(c_{q,p+1-q}) . $$
\end{enumerate}
Then, the proof goes as follows. Let us consider
\begin{eqnarray}\nonumber
    u^{(p-q)} = u^{(p+1-q)} +
    \frac{(-1)^{p+2}(-1)^{p+1}\ldots(-1)^{p+2-q}}{(p+1)p(p-1)\ldots(p+2-(q+1))}W(c_{q,p-q})
    .
\end{eqnarray}
Then, by $(c')_q$, by Lemma \ref{W2} and by $(a')_{q+1}$, we
obtain $\|u^{(p-q)}\|^2 \leq c'_3 \|\delta c\|^2$. Moreover, by
$(c')_q$ and Lemma \ref{W1} we have
\begin{eqnarray}\nonumber
    \frac{(-1)^{p+2}}{p+1} d W(c) &=&
    d u^{(p-q)} - \frac{(-1)^{p+2}\ldots(-1)^{p+2-q}}{(p+1)\ldots(p+2-(q+1))}
    d W(c_{q,p-q})  \\ \nonumber &&+ \frac{(-1)^{p+2}\ldots(-1)^{p+2-q}}{(p+1)\ldots(p+2-q)}
    W(c_{q,p+1-q})\\ \nonumber
    &=&
    d u^{(p-q)} - \frac{(-1)^{p+2}\ldots(-1)^{p+2-q}}{(p+1)\ldots(p+2-q)}
    W(c_{q,p+1-q})\\ \nonumber && + \frac{(-1)^{p+2}\ldots(-1)^{p+2-q}}{(p+1)\ldots(p+2-(q+1))}
    (-1)^{p+1-q}W(c_{q+1,p-q}) \\ \nonumber &&
    + \frac{(-1)^{p+2}\ldots(-1)^{p+2-q}}{(p+1)\ldots(p+2-q)}
    W(c_{q,p+1-q})\nonumber
\end{eqnarray}
and the claim follows.

\smallskip

At the end of the induction (for $q=p+1$), we obtain a $p$-form
$u^{(0)}$ such that $\|u^{(0)}\|^2 \leq c'_3 \|\delta c\|^2$ and
\begin{eqnarray} \nonumber
    d W(c) & =&
    (-1)^p(p+1) \left(d u^{(0)} + k(p) W(c_{p+1,0})\right)\\ \nonumber
    & =&
    (-1)^p(p+1) \left(d u^{(0)} + k(p) W(r(d \S(c)))\right)\\ \nonumber
    & =&
    (-1)^p(p+1) \left(d u^{(0)} + k(p) d(\S(c))\right)
\end{eqnarray}
where $k(p)$ is a constant depending only on $p$. Moreover, as
$\S(c)$ is a coexact $p$-form, this implies
\begin{eqnarray}\nonumber
    coex(W(c)) = (-1)^p(p+1) \left(coex( u^{(0)} )+ k(p) \S(c)\right)
\end{eqnarray}
so that
\begin{eqnarray}\nonumber
    \|coex(W(c))\|  &\geq& (p+1)|k(p)|\cdot \|\S(c)\| - (p+1)\|u^{(0)}\|\\
                    &\geq& (p+1)|k(p)|\cdot \|\S(c)\| -
                    (p+1)(c'_3)^{\frac{1}{2}}\|\delta c\| .
                    \nonumber
\end{eqnarray}
But, by $(ii)_{\S}$, if $\|\delta c\|^2 \leq \Lambda' \|c\|^2$
then $\|\S(c)\| \geq (c'_2)^{\frac{1}{2}}\|c\|$. Therefore,
\begin{eqnarray}\nonumber
    \|coex(W(c))\|  &\geq& (p+1)|k(p)|(c'_2)^{\frac{1}{2}}
    \left( \|c\| - \sqrt{\frac{c'_3}{c'_2 k(p)^2}}\|\delta c\|\right)
\end{eqnarray}
Finally, if $\|\delta c \|^2 \leq \tilde{\Lambda} \|c\|^2$, with
$\tilde{\Lambda} = \min\left\{\Lambda',
\frac{k(p)^2c'_2}{4c'_3}\right\}$, then
\begin{eqnarray}\nonumber
    \|coex(W(c))\|  \geq
    \frac{1}{2}(p+1)|k(p)|(c'_2)^{\frac{1}{2}}
     \|c\|
\end{eqnarray}
which is the desired inequality in $(ii)_{\tilde{S}}$. $\square$


\subsection{Another proof of "McGowan lemma"} \label{sec43}

In \cite{MG}, the author gives a lower bound for the
$N^{\text{th}}$ eigenvalue of $\Delta$ on exact 2-forms on a
compact Riemannian manifold $M$ (see Lemma 2.3 in \cite{MG}) where
$N$ depends on an open cover of $M$. In particular, if the open
cover is contractible then $N-1$ is the number of non-empty
intersections of triples of open sets in the open cover. The lower
bound depends then essentially on lower bounds for the smallest
positive eigenvalue of $\Delta$ on exact forms on the open sets of
the cover, on the intersection of pairs of such open sets and on
the intersection of triples of such open sets. The proof of
McGowan relies also on the double complexe of \Cech- de Rham and
can be compared to the induction done in Section \ref{constrD} to
construct the discretizing operator $\D$. So it is not so
surprising that we obtain the following generalization of the
lemma. The main difference is that in our technique, if the
discretization is of sufficiently small mesh then Lemma
\ref{lmmmu} gives the lower bound for the spectrum on the
intersections. But, then $N$ can get quite large as it is
comparable to the number of open sets in the open cover. Let us
now state and prove the result.

\begin{lmm}\label{lmmmcg}
    Let $n\geq 1$, $\kappa \geq 0$ and $r_0 > 0$. Then there
    exists a positive constant $\lambda(n,\kappa,r_0)$ depending
    only on $n$, $\kappa$ and $r_0$ such that for any $(M,g)\in
    \K$ we have $$\lambda_{N,p}^{d^*}(M) \geq
    \lambda(n,\kappa,r_0)$$ where $N \leq c(n,p,\kappa,r_0) Vol(M)$ and
    $c(n,p,\kappa,r_0)$ is a positive constant.
\end{lmm}

\textbf{Proof}: let $\rho_0$ be given by Lemma \ref{lmmrho} and
let $X$ a $\frac{\rho_0}{4}$-discretization of $M$. Then the
discretizing operator $$\D : d^* \Lambda^{p+1}(M) \to
\delta^*\CqR{p+1}{\U}$$ constructed in Section \ref{constrD}
satisfies $(i)_{\D}$ and $(ii)_{\D}$. Let then $$N = dim \left(
\delta^*\CqR{p+1}{\U} \right) + 1. $$ Consider moreover $\phi_1,
\ldots, \phi_N$ the $N$ first eigenforms in $d^*
\Lambda^{p+1}(M)$. By definition of $N$, there exist $a_1, \ldots,
a_N$ such that $\sum_{i=1}^N a_i \D(\phi_i)=0$ and $\sum_{i=1}^N
a_i \phi_i \neq 0$. In particular, by $(ii)_{\D}$, we get
$$\left\| d\left(\sum_{i=1}^N a_i \phi_i \right)\right\|^2 \geq
\Lambda \left\| \sum_{i=1}^N a_i \phi_i \right\|^2$$ and thanks to
the variational characterization of the spectrum
$$\lambda_{N,p}^{d^*}(M) =
\max\left\{\frac{\|d\phi\|^2}{\|\phi\|^2} : \phi \in \langle
\phi_1, \ldots, \phi_N \rangle \setminus \{0\} \right\} \geq
\Lambda .$$ Note that by Remark \ref{cardX}, we have $N \leq
|\S_p(\U_X)| \leq  c_2 \frac{\nu^p}{(p+1)!} Vol(M)$ where $c_2$
and $\nu$ depend only on $n$, $p$, $\kappa$ and $r_0$. $\square$

%% file: DiscFormsAppendix.tex
\appendix

\section{Appendix}

\subsection{Whitney forms} \label{whitney}

Let $(M^n,g)$ be a compact connected $n$-dimensional Riemannian
manifold without boundary. Let $\U$ be a finite contractible open
cover of $M$. Let $\{\varphi_j\}$ be a partition of unity
subordinated to $\U$. Let $\nu$ a bound on the cardinality of the
sets $\{j : U_j \cap U_I \neq \emptyset\}$, $I\in S_q(\U)$, $q
\geq 0$.

\begin{dfn}
For any $I = \{i_0,\ldots,i_q\} \in S_q(\U)$, we define the
\textbf{Whitney form} $W_I \in \Lambda^q(M)$ by $$W_I =
\sum_{j=0}^{q} (-1)^j \varphi_{i_j} d \varphi_{i_0} \wedge \ldots
\wedge d \varphi_{i_{j-1}} \wedge d \varphi_{i_{j+1}} \wedge
\ldots \wedge d \varphi_{i_q}$$
\end{dfn}

\begin{rmk}
Note that $W_I$ has support in $U_I$. Moreover, we have $dW_I =
(q+1) d\varphi_{i_0} \wedge \ldots \wedge d \varphi_{i_q}$, for $I
= \{i_0,\ldots,i_q\}$. In the sequel, we will write $d\varphi_I =
d\varphi_{i_0} \wedge \ldots \wedge d \varphi_{i_q}$.
\end{rmk}

We can extend the definition of Whitney forms to $q$-cochains as
follows.

\begin{dfn} Let $ W : \Cq{q}{U}{\Lambda^p} \to \Lambda^{p+q}(M)$
the application defined by $$W(c) = \sum_{I \in S_q(\U)} W_I
\wedge c(I) . $$
\end{dfn}
The application $W$ restricted to \Cech cochains is the Whitney
map introduced by Whitney (see \cite{Wh}) (up to a constant). The
following lemma generalizes the well-known fact that the Whitney
map commutes with the exterior differential and the coboundary.

\begin{lmm}\label{W1} For any $c \in \Cq{q}{U}{\Lambda^p}$, we have $$d
W(c)= (q+1) W (\delta c) + (-1)^q W (d c).$$
\end{lmm}
\textbf{Proof}: we have
\begin{eqnarray}
    d W(c)  &=& \sum_{I \in S_q(\U)} d (W_I \wedge c(I)) \nonumber\\
            &=& \sum_{I \in S_q(\U)} d W_I \wedge c(I) +
            (-1)^q\sum_{I \in S_q(\U)} W_I \wedge d c(I) \nonumber\\
            &=& (q+1) \sum_{I \in S_q(\U)} d \varphi_I \wedge c(I)
            + (-1)^q W(d c) . \nonumber
\end{eqnarray}
Let us now compute $W(\delta c)$. We have
\begin{eqnarray}\nonumber
    W(\delta c) &=& \sum_{J \in S_{q+1}(\U)}
                W_J \wedge \left(\sum_{j\in J}\epsilon(j,J\setminus j) c(J\setminus j)\right)
\end{eqnarray}
where $\epsilon(j,J\setminus j)$ is $\pm 1$ according to the
signature of the permutation ordering the set $\{j\} \cup
(J\setminus j)$ in $J$. If we let $I = J \setminus j$, we can
write
\begin{eqnarray} \nonumber
   W(\delta c) &=& \sum_{I \in S_{q}(\U)} \sum_{j : U_j \cap U_I
   \neq \emptyset}  W_{\{j,I\}} \wedge  c(I)
\end{eqnarray}
so that it suffices to show that
\begin{eqnarray} \label{fin}
\sum_{j : U_j \cap U_I \neq \emptyset}  W_{\{j,I\}} =d\varphi_I
\end{eqnarray}
to conclude the proof. Let us rewrite the expression as follows
\begin{eqnarray}\label{fin2}
\sum_{j : U_j \cap U_I \neq \emptyset}  W_{\{j,I\}}
    = \sum_{j : U_j \cap U_I \neq \emptyset} \varphi_j d\varphi_I
    - d\varphi_j \wedge W_I .
\end{eqnarray}
But as $\{\varphi_j\}$ is a partition of unity $  \sum\limits_{j :
U_j \cap U_I \neq \emptyset} \varphi_j =1$ and $ \sum\limits_{j :
U_j \cap U_I \neq \emptyset} d\varphi_j =0$, hence (\ref{fin2})
implies (\ref{fin}). $\square$

\begin{lmm}\label{W2} There exists a positive constant $k$
depending only on $n$, $\nu$ and on $\|d\varphi_j\|_{\infty}$ such
that for any \Cech cochain $c$, $\| W(c)\|^2 \leq k \|c\|^2$.
\end{lmm}
\textbf{Proof}: it follows from the definition of $W$ and from a
direct calculation. $\square$


\subsection{About the convexity of balls}\label{pflmmrho}

\textbf{Proof of Lemma \ref{lmmrho}}: the main idea to prove this
lemma is to smooth $g$ to obtain a more regular metric $\tilde{g}$
and then compare $\tilde{g}$ to a Euclidean metric $\tilde{e}$. We
do not compare directly $g$ with a Euclidean metric as we need to
control the difference between the different connections involved.
So let $(M,g) \in \K$. It follows from a result of Abresch (see
\cite{CFG}, Theorem 1.12) that there exists a Riemannian metric
$\tilde{g}$ on $M$ such that
\begin{enumerate}
\item[(a)] $e^{-\frac{1}{4}} g \leq \tilde{g} \leq e^{\frac{1}{4}} g$
\item[(b)] $|\nabla^g - \nabla^{\tilde{g}}| \leq \frac{1}{4}$
\item[(c)] $|K_{\tilde{g}}| \leq \tilde{\kappa}(n,\kappa) \text{
and }|\nabla^{\tilde{g}} R_{\tilde{g}}| \leq k(n,\kappa)$
\end{enumerate}
where $\tilde{\kappa}$ and $k$ depend only on $n$ and $\kappa$. In
particular, (a) implies that, the length of the curves, the
distances and the volumes are comparable within a ratio depending
only on $n$. Moreover, if $B$ denotes a ball for $g$ and
$\tilde{B}$ a ball for $\tilde{g}$, we get $B(p,e^{-\frac{1}{2}}r)
\subseteq \tilde{B}(p,r) \subseteq B(p,e^{\frac{1}{2}}r)$. First,
we show that there exists $\tilde{r}_0 >0$ depending only on $n$,
$\kappa$, $r_0$ such that
\begin{eqnarray}\label{inj}
 inj(M,\tilde{g}) \geq \tilde{r}_0 .
\end{eqnarray}
This is a direct consequence of a theorem of Klingenberg and a
theorem of Cheeger. Indeed, by Klingenberg's Theorem (see for
instance \cite{Be}, Theorem 89, or \cite{Ko} and \cite{Ka}) and as
we have bounded sectional curvature, the injectivity radius
satisfies
\begin{eqnarray} \label{kling}
inj(M,\tilde{g}) \geq \min \left\{
\frac{\pi}{\sqrt{\tilde{\kappa}}} , \frac{1}{2}
\tilde{l}(\tilde{\gamma}) \right\}
\end{eqnarray}
where $\tilde{l}(\tilde{\gamma})$ is the length with respect to
$\tilde{g}$ of the shorter smooth geodesic (w.r.t. $\tilde{g}$)
loop. Moreover, we show there exists $L > 0$ depending only on
$n$, $\kappa$ and $r_0$ such that $\tilde{l}(\tilde{\gamma}) \geq
L $ as follows. First, if $\frac{1}{2}\tilde{l}(\tilde{\gamma})
\geq e^{-\frac{1}{2}} r_0$,  $L = 2 e^{-\frac{1}{2}} r_0$ is
suitable. Then suppose $\frac{1}{2}\tilde{l}(\tilde{\gamma}) <
e^{-\frac{1}{2}} r_0$. By construction of $\tilde{g}$ we get then
$\frac{1}{2}l(\tilde{\gamma}) <  r_0$ i.e. $\tilde{\gamma}$ is
contained in $B(\tilde{\gamma}(0),r_0) = B$. Again by construction
of $\tilde{g}$, $\widetilde{Vol}(B) \geq c(n) Vol(B)$ and as
$(M,g) \in \K$, there exists $c(n,\kappa, r_0)>0$ such that
$Vol(B) \geq c(n,\kappa,r_0)$ so that $\widetilde{Vol}(B)$ is
bounded below by a constant $V$ depending only on $n$, $\kappa$
and $r_0$. Moreover, $\widetilde{diam}(B) \leq 2
e^{\frac{1}{2}}r_0 = d$. So that we can apply Theorem 2.1. of
Cheeger in \cite{Che} that ensures the existence of a positive
constant $L$ depending only on $d$, $V$ and $\tilde{\kappa}$ and
therefore only on $n$, $\kappa$ and $r_0$ such that
$\tilde{l}(\tilde{\gamma}) \geq L$. Together with ($\ref{kling}$),
this implies (\ref{inj}).

\smallskip

A suitable candidate to be the diffeomorphism cited in the claim
is the exponential map with respect to the metric $\tilde{g}$. Let
then $$\varphi = \widetilde{\exp}_p : B(0,\tilde{r}_0) \to
\tilde{B}(p,\tilde{r}_0)$$ and $\tilde{e}$ the Euclidean metric on
$\tilde{B}(p,\tilde{r}_0)$ induced by $\varphi^{-1}$ and the
normal coordinates. As soon as $e^{\frac{1}{2}}r \leq
\tilde{r}_0$, we have $B(p,r)\subseteq \tilde{B}(p,\tilde{r}_0)$
and then $\varphi^{-1}(B(p,r))$ is well-defined. We are going to
show now that there exists a positive constant $0 <
\rho_0(n,\kappa,r_0) \leq e^{-\frac{1}{2}}\tilde{r}_0$ such that
for any $B(q,\rho) \subseteq B(p,\rho_0) \subseteq
\tilde{B}(p,\tilde{r}_0)$ we have
\begin{eqnarray}
\varphi^{-1}(B(q,\rho)) \text{ is a Euclidean convex.}
\end{eqnarray}
This is equivalent to showing that the application
\begin{eqnarray}
f : (B(q,\rho),\tilde{e}) \to \R \text{ , } x \mapsto
\frac{1}{2}d(q,x)^2
\end{eqnarray}
is convex (w.r.t. $\tilde{e}$), in other words that the Hessian of
$f$ with respect to $\tilde{e}$ is non-negative i.e. $
D^2_{\tilde{e}}f(U,U) \geq 0$ on $B(q,\rho)$, for $\rho$ and
$\rho_0$ well-chosen. Let us recall the following definition of
the Hessian $$D^2f(U,V) = U \cdot df(V) - df(\nabla_U V)$$ where
$\nabla$ is the Levi-Civita connection. Using this definition of
the Hessian for $\tilde{e}$ and $g$, we get
\begin{eqnarray}\label{hess}
D^2_{\tilde{e}}f(U,U)& = & D^2_{g}f(U,U) + df (\nabla^g_U U -
\nabla^{\tilde{e}}_U U  ) \nonumber \\ &=& D^2_{g}f(U,U) + df
(\nabla^g_U U - \nabla^{\tilde{g}}_U U) + df (\nabla^{\tilde{g}}_U
U - \nabla^{\tilde{e}}_U U ) .
\end{eqnarray}
Proposition 6.4.6. of Buser and Karcher in \cite{BK} says that
\begin{eqnarray}\nonumber
D^2_{g}f(U,U) \geq \rho \frac{s'_{\kappa}(\rho)}{s_{\kappa}(\rho)}
g(U,U)
\end{eqnarray}
where $s_{\kappa}(\rho) = \frac{1}{\sqrt{\kappa}}
\sin(\sqrt{\kappa} \rho)$. So that
$\frac{s'_{\kappa}(\rho)}{s_{\kappa}(\rho)} =
\sqrt{\kappa}\cot(\sqrt{\kappa}\rho)$ and hence there exists
$\rho_1(\kappa) > 0$ such that for any $0< \rho < \rho_1$,
$\frac{s'_{\kappa}(\rho)}{s_{\kappa}(\rho)} \geq 1$. Therefore, on
$B(q,\rho)$ with $\rho \leq \rho_1$ we have
\begin{eqnarray}\label{hessbk}
D^2_{g}f(U,U) \geq \rho g(U,U)
\end{eqnarray}
and this shows also that for such $\rho$'s, $B(q,\rho)$ is convex
(w.r.t. $g$). Also as a consequence of Proposition 6.4.6. of
\cite{BK}, we get
\begin{eqnarray}\label{gradbk}
g(\nabla^g f, \nabla^g f) \leq \rho^2
\end{eqnarray}
where $\nabla^g f$ is the gradient of $f$ with respect to $g$.

\smallskip

Moreover, by construction of $\tilde{g}$ and by $(b)$ in the
result of Abresch, we have
\begin{eqnarray}\label{nablab}
|\nabla^g_U U - \nabla^{\tilde{g}}_U U|_g \leq \frac{1}{4} g(U,U)
\text{.}
\end{eqnarray}

\smallskip

By construction of $\tilde{e}$ and as the
$\nabla^{\tilde{g}}R_{\tilde{g}}$ is uniformly bounded, Corollary
1 of Kaul in \cite{Kau} asserts the existence of an application $h
\geq 0$ such that
\begin{eqnarray}\nonumber
|\nabla^{\tilde{g}}_U U - \nabla^{\tilde{e}}_U U|_{\tilde{g}} (y)
\leq h(\tilde{d}(p,y)) \tilde{g}(U,U)
\end{eqnarray}
with $h(0) = 0$ and $h$ depends only on bounds on $K_{\tilde{g}}$
and $\nabla^{\tilde{g}}R_{\tilde{g}}$. Hence, there exists
$R(n,\kappa,r_0) > 0$ such that for any $r \leq R$, $h(r) \leq
\frac{1}{4}e^{-\frac{3}{4}}$. So that we obtain on
$\tilde{B}(p,r)$ with $r \leq R$
\begin{eqnarray}\label{nablka}
|\nabla^{\tilde{g}}_U U - \nabla^{\tilde{e}}_U U|_g \leq
e^{\frac{1}{2}} |\nabla^{\tilde{g}}_U U - \nabla^{\tilde{e}}_U
U|_{\tilde{g}} \leq \frac{1}{4} e^{-\frac{1}{4}} \tilde{g}(U,U)
\leq \frac{1}{4}g(U,U) \text{.}
\end{eqnarray}

\smallskip

Finally, introduce (\ref{hessbk}), (\ref{gradbk}), (\ref{nablab})
and (\ref{nablka}) in (\ref{hess}) and let us define $\rho_0 =
\min \{e^{-\frac{1}{2}}\tilde{r}_0 , \rho_1 , e^{-\frac{1}{2}}R\}$
to obtain the following. We have $B(p,\rho_0)\subseteq
\tilde{B}(p,\tilde{r}_0)$, $B(p,\rho_0)\subseteq \tilde{B}(p,R)$
and for any $B(q,\rho) \subseteq B(p,\rho_0)$, $\rho \leq \rho_1$
holds. Hence on $B(p,\rho_0)$ and for any $B(q,\rho) \subseteq
B(p,\rho_0)$ we have
\begin{eqnarray}
D^2_{\tilde{e}}f(U,U) \geq \rho g(U,U) - \frac{1}{4}\rho g(U,U) -
\frac{1}{4}\rho g(U,U) = \frac{1}{2}\rho g(U,U) \geq 0
\end{eqnarray}
i.e. $f$ is convex. To conclude the proof, we remark that
$$B(q,\rho) \subseteq B(p,d(p,q)+\rho) \subseteq
\tilde{B}(p,e^{\frac{1}{2}}(d(p,q)+\rho)) \subseteq
\tilde{B}(p,\tilde{r}_0)$$ so that $\varphi^{-1}$ restricted to
$\tilde{B}(p,e^{\frac{1}{2}}(d(p,q)+\rho))$ is a quasi-isometry
with constants of quasi-isometry depending only on $n$, $\kappa$
and $d(p,q)+ \rho$. More precisely,
$(\tilde{B}(p,e^{\frac{1}{2}}(d(p,q)+\rho)),\tilde{e})$ is
quasi-isometric to
$(\tilde{B}(p,e^{\frac{1}{2}}(d(p,q)+\rho)),\tilde{g})$ with
constants of quasi-isometry depending only on $d(p,q)+ \rho$ and
$\tilde{\kappa}(n,\kappa)$ and by construction of $\tilde{g}$ we
can deduce that $(B(q,\rho),g)$ is quasi-isometric to
$(B(q,\rho),\tilde{e})$ with constants of quasi-isometry depending
only on $n$, $\kappa$ and $d(p,q)+ \rho$. This ends the proof of
the lemma. $\square$


%% file: DiscFormsReferences.tex
\smallskip

{\small Tatiana Mantuano

Universit\'e de Neuch\^atel

Institut de Math\'ematiques

rue Emile-Argand 11

2009 Neuch\^atel

Switzerland

e-mail: Tatiana.Mantuano@unine.ch }